\documentclass[12pt]{article}
\input epsf.tex
\usepackage{latexsym}
\usepackage{amsmath,amsthm,amsfonts,amssymb}
\usepackage{epsfig}

\newtheorem{pro}{Proposition}[section]
 \newtheorem{thm}[pro]{Theorem}
 \newtheorem{lem}[pro]{Lemma}
 \newtheorem{defn}[pro]{Definition}

 \newtheorem{cnj}[pro]{Conjecture}
 \newtheorem{cor}[pro]{Corollary}
 \newtheorem{question}[pro]{Question}
 \newtheorem{remark}[pro]{Remark}

 \def\B{{\cal B}}
 \def\C{{\mathbb C}}

 \def\H{{\mathbb H}}
 \def\G{{\cal G}}

 \def\P{{\cal P}}

 \def\M{{\cal M}}
\def\m{{\bar M}}

 \def\N{{\mathbb N}}
 \def\R{{\mathbb R}}
 \def\chix{{\raise.5ex\hbox{$\chi$}}}

\def\mE{{\mathbb E}}

\def\bgamma{{\bar \gamma}}

\def\sO{{\cal O}}
\def\sF{{\cal F}}

\def\sH{{\cal H}}

\def\sB{{\cal B}}

\def\sK{{\cal K}}

\def\sV{{\cal V}}
\def\sU{{\cal U}}
\def\tH{{\tilde H}}

\def\tM{{\tilde M}}

\def\tgamma{{\tilde \gamma}}

\def\tp{{\tilde p}}
\def\tq{{\tilde q}}
\def\tv{{\tilde v}}
\def\tw{{\tilde w}}
\def\hH{{\hat H}}

\def\wH{{H}}

\def\arccosh{{\textnormal{arccosh}}}
\def\gM{{M}}
\def\m{{\hat M}}

\begin{document}
\title{Immersions of Pants into a Fixed Hyperbolic Surface}
\author{Lewis Bowen}
\maketitle
\begin{abstract}
%Let $S$ be a fixed closed hyperbolic surface. The main objects of
%study in this paper are locally isometric immersions $j:P \to S$ from
%pairs of pants $P$ into $S$. Asymptotic estimates for the number of
%such immersions satisfying various constraints are obtained. We also show that, in an asymptotic sense, given two
%closed geodesics $\gamma_1, \gamma_2$ of length close to $L$, there
%is almost always an immersed $4$-holed sphere
%partially bounded by $\gamma_1 \cup \gamma_2$ whose other boundary
%components also have length close to $L$. The methods generalize to immersions of other surfaces with boundary into closed hyperbolic $n$-manifolds for $n=2,3$.%

Exploiting a relationship between closed geodesics on a generic closed hyperbolic surface $S$ and a certain unipotent flow on the product space $T_1(S) \times T_1(S)$, we obtain a local asymptotic equidistribution result for long closed geodesics on $S$. Applications include asymptotic estimates for the number of pants immersions into $S$ satisfying various geometric constraints. Also we show that two
closed geodesics $\gamma_1, \gamma_2$ of length close to $L$ chosen uniformly at random have a high probability of partially bounding an immersed $4$-holed sphere whose other boundary
components also have length close to $L$.

\end{abstract}

\section{Introduction}

One motivation for this work came from two well-known conjectures:

\begin{cnj}(The Surface Subgroup Conjecture)
Let $\M$ be a closed hyperbolic $3$-manifold. Then there exists a $\pi_1$-injective map $j:S\to \M$ from a closed surface $S$ of genus at least 2 into $\M$.
\end{cnj}

\begin{cnj}(The hyperbolic Ehrenpreis conjecture)(\cite{Ehr},\cite{Gen})
Let $\epsilon >0$ and let $S_1, S_2$ be two closed hyperbolic surfaces. Then there exists finite-sheeted locally isometric covers  ${\tilde S_i}$ of $S_i$ (for $i=1,2$) such that there is a $(1+\epsilon)$ bi-Lipschitz homeomorphism between ${\tilde S_1}$ and ${\tilde S_2}$.
\end{cnj}

%A bolder conjecture than the first above is that $j$ can, in addition, be chosen to be ``almost geodesic'' (meaning, for example, that the principal curvatures of the image surface are all bounded close to 0). This bolder statement is directly analogous to the second conjecture above: if there exists such finite covers with a $(1+\epsilon)$ bi-Lipshitz homeomorphism $h:{\tilde S_1} \to {\tilde S_2}$ between them then the graph of $h$ in ${\tilde S_1} \times {\tilde S_2}$ pushes down to an immersed surface in the product $S_1 \times S_2$ that is locally ``almost conjugate into the diagonal''. A converse also holds although we will make neither statement precise here.

Attempting to understand these conjectures led to the study of immersions of three-holed spheres into $3$-manifolds $\M$ and cross products $S_1 \times S_2$; we would like to glue such immersions together to obtain an immersion of a closed surface $S \to \M$ or $S\to S_1 \times S_2$ with tight control on its geometric structure.

But the simpler case, that of immersions of three-holed spheres into a closed hyperbolic surface $S$, is not well-understood. For instance, the following question is unknown.

\begin{question}\label{cnj:singlesurface}
Given $\epsilon>0$, for sufficiently large $L$, does there exist a finite sheeted cover $\pi: {\tilde S} \to S$ such that ${\tilde S}$ admits a pair of pants decomposition, every geodesic of which has length in $(L-\epsilon,L+\epsilon)$?
\end{question}

Another motivation for the present work comes from a desire to ``bridge the gap'' between two research areas: group growth and subgroup growth. The former concerns itself with the asymptotic number of elements in a given group of word length less than $R$, the latter with the asymptotic number of subgroups with finite index less than $R$. Is there something in between? Among other things, here we study the asymptotic number of conjugacy classes of 2-generator subgroups of a surface group satisfying certain geometric conditions with the aim (not yet realized) of amalgamating these subgroups together to obtain finite index subgroups with geometric constraints.

%Another motivation is to shed some light on the ``higher-dimensional'' spectra of a hyperbolic surface $S$. To be precise, let ${\tilde \Sigma}$ be the set of all locally-isometric immersions $i=(j: P \to S)$ from a genus $g$, $n$-holed hyperbolic surface $P$ into $S$. For each $(j:P\to S)\in {\tilde \Sigma}$ let $mod(i)$ be the modulus of $P$ in the moduli space $M(g,n)$ of all genus $g$, $n$-holed surfaces. The set $\Sigma(g,n)=mod({\tilde \Sigma})$ is the $(g,n)$-spectrum of $S$. For example, the $(0,2)$-spectrum is the length spectrum. A fundamental issue taken up here is the asymptotic growth rate of the $(0,3)$-spectrum. But more general questions deserve to be asked: what is the relationship between the various spectra? It is obvious that the $(0,3)$ spectrum determines the $(0,2)$-spectrum. Is the converse true?  %

%Our main purpose is to demonstrate techniques for studying immersions of 3-holed spheres into a fixed hyperbolic surface. These methods can be generalized to the case of $3$-manifolds, products $S_1 \times S_2$ and immersions of other surfaces with boundary but we do not take that up here. 

\subsection{Equidistibution Results}

Let $S$ be a fixed closed hyperbolic surface. Regard $S$ as the quotient space $\H^2/\Gamma$ where $\Gamma$ is a fixed lattice in the group $Isom^+(\H^2)$ $(=PSL_2(\R))$ of all orientation-preserving isometries of the hyperbolic plane $\H^2$.

Let $w$ be an arbitrary unit vector in the unit tangle bundle $T_1(S)$. Consider the geodesic segment of length $L$ tangent to $w$ with $w$ based at its midpoint. If the tangent vectors $e_1, e_2$ at the endpoints are close then a short segment can be adjoined to it to obtain a closed path in $S$. The closed geodesic $\gamma$ in the homotopy class of this path is very close to the original segment. A calculation we will use often quantifies how close. For example, we show that there is a function $F=(F_1,F_2,F_3)$ of the position of $e_1$ relative to $e_2$ such that the distance from $w$ to $\gamma$ along a geodesic segment orthogonal to $w$ equals $F_1 e^{-L/2} + O(e^{-L})$, the angle at which this segment intersects $\gamma$ equals $\pi/2 + F_2 e^{-L/2} + O(e^{-L})$ and the length of $\gamma$ equals $L+ F_3 + O(e^{-L})$. Roughly speaking, if the distance between $e_1$ and $e_2$ is less than $\epsilon$ then $|F_1|,|F_2|,|F_3|$ are all less than $\epsilon$ as well. For a precise statement see corollaries \ref{cor:perp} and \ref{cor:2d} below.

If the vectors $e_1$ and $e_2$ are not close then we push $w$ to its right along an orthogonal geodesic. The pair $(e_1,e_2)$ moves by a product of hypercycle flows in the product space $T_1(S) \times T_1(S)$. Pushing $w$ a distance $O(e^{-L/2})$ amounts to flowing this pair for $O(1)$ time. As $L$ tends to infinity, this product hypercycle flow converges to a product of horocycle flows. As a consequence of Ratner's work on Raghunathan's conjectures we show that if the commensurator, $Comm(\Gamma)$, contains only orientation preserving isometries then the latter flow is uniformly equidistributed on the product space. By definition, $Comm(\Gamma)$ is the set of all isometries $g \in Isom(\H^2)$ such that $g\Gamma g^{-1} \cap \Gamma$ has finite index in $\Gamma$.

To make this precise, for $w\in T_1(\H^2)$, $T,L>0$, let $\mu = \mu_{w,T,L}$ be the probability measure on the set 
\begin{eqnarray*}
\big\{(w_t G_{-L/2}, w_t G_{L/2}) \in T_1(\H^2)\times T_1(\H^2)\, \big| \, 0 \le t \le T e^{-L/2}\big\}
\end{eqnarray*}
induced by Lebesgue measure on $[0,Te^{-L/2}]$. Here $w_t$ is the unit vector obtained by pushing $w$ to its right along an orthogonal geodesic for time $t$ at unit speed. $G_L$ is the geodesic flow for time $L$. So if $\sigma_t$ is the segment of length $L$ tangent to $w_t$ with $w_t$ based at its midpoint, then $w_t G_{-L/2}$ and $w_t G_{L/2}$ are the unit vectors tangent to $\sigma_t$ at its ends and oriented consistently with $w_t$. Let $\pi_*(\mu_{w,T,L})$ be the projection of $\mu_{w,T,L}$ to $T_1(S) \times T_1(S)$.

\begin{thm}\label{thm:equidistribution}
Assume $Comm(\Gamma)<Isom^+(\H^2)$. Let $\lambda$ be Haar probability measure on $T_1(S)$. Given any continuous function $f: T_1(S) \times T_1(S) \to \C$ and any $\epsilon_0 > 0$ there exists a $T_0 > 0$ such that for all $T, L>T_0$ and $w \in T_1(\H^2)$ 
\begin{eqnarray*}
|\pi_*(\mu_{w,T,L})(f)-\lambda \times \lambda(f)| < \epsilon_0.
\end{eqnarray*}
 
 \end{thm}
%This translates into a local asymptotic equidistribution result for long closed geodesics on $S$. We use this to build pants immersions into $S$ by constructing generators for the image subgroup corresponding to two of the boundary components.

 The condition $Comm(\Gamma)<Isom^+(\H^2)$ is generic in the moduli space of genus $g$ surfaces for every genus $g$ (see remark \ref{rmk:margulis}). With the calculation results, this theorem implies local asymptotic equidistribution for long closed geodesics on $S$. 

To give a sample of what can be obtained, for $l_1< l_2$ let $\G_L(l_1,l_2)$ be the set of all closed oriented geodesics on $S$ with length in the interval $(L+l_1,L+l_2)$. We allow closed geodesics to cover their images multiple times. So the length of a geodesic $\gamma$ equals $m$ times the length of its image $\bgamma$ for some $m \in \N$.  The next result concerns the asymptotic number of almost-perpendicular intersections of the union of geodesics in $\G_L(l_1,l_2)$ with a sequence of segments $\sigma_L$ on the surface.

%The commensurator of $\Gamma$, $Comm(\Gamma)$ is the set of all isometries $g\in Isom(\H^2)$ such that $g\Gamma g^{-1} \cap \Gamma$ has finite index in $\Gamma$.

\begin{thm}\label{thm:angle}
Assume $Comm(\Gamma)<Isom^+(\H^2)$. Let $\{\sigma_L\}_{L >0}$ be a sequence of oriented geodesic segments $\sigma_L \subset S$. Assume $length(\sigma_L)e^{L/2} \to \infty$ as $L\to \infty$. Let $(l_1,l_1), (a_1,a_2)$ be finite intervals of the real line. Let $N_L= N_L(a_1,a_2,l_1,l_2)$ be the set of unit vectors $v$ such that
\begin{itemize}
\item the basepoint of $v$ is in $\sigma_L$,
\item $v$ is tangent to a closed geodesic $\gamma \in \G_L(l_1,l_2)$,
\item $v$ is oriented consistently with $\gamma$,
\item the angle from $v$ to $\sigma_L$ is in the interval $\pi/2 + (a_1,a_2)e^{-L/2}$.
\end{itemize}
Then 
\begin{eqnarray*}
\# N_L \sim \frac{length(\sigma_L)}{vol(T_1(S))}(a_2-a_1)(e^{l_2}-e^{l_1})e^{L/2}.
\end{eqnarray*}
\end{thm}
Here and throughout the paper, $F\sim G$ means $\lim_{L \to \infty} \frac{F}{G} = 1$. $vol(T_1(S))$ equals $2\pi area(S) = (2\pi)^2(2genus(S)-2)$. We also give a new proof of a special case of Rufus Bowen's equidistribution theorem (see theorem \ref{thm:bowen}).

\subsection{Counting Pants Immersions}
We use the theorems above to build and count pants immersions into $S$ by constructing generators for the image subgroup corresponding to two of the boundary components. To state the results, fix $\epsilon>0$. For $r_1,r_2,r_3,L>0$ let $\P_L(r_1,r_2,r_3)$ be the set of all locally isometric orientation-preserving immersions $j:P \to S$ in which $P$ is a hyperbolic three-holed sphere (i.e. a pair of pants) with geodesic boundary components of length $l_1,l_2,l_3$ satisfying 
\begin{eqnarray*}
l_i \in r_iL + (-\epsilon,\epsilon)
\end{eqnarray*}
for $i=1,2,3$. We implicitly identify immersions $j_1:P_1 \to S$ and $j_2:P_2 \to S$ if there is an isometry $\Psi: P_1 \to P_2$ such that $j_1 = j_2 \circ \Psi$. Thus $\P_L(r_1,r_2,r_3)$ is a finite set. The first result gives asymptotics for the cardinality $\P_L(r_1,r_2,r_3)$.

%If $f,g$ are functions of $\epsilon$, we write $f \approx g$ to if there exists positive universal constants $k_1,k_2$ such that
%\begin{eqnarray*}
%k_1 F(\epsilon) \le G(\epsilon) \le k_2 F(\epsilon).
%\end{eqnarray*}

\begin{cor}\label{cor:all}
 Assume $Comm(\Gamma) < Isom^+(\H^2)$. Then
\begin{eqnarray*}
|\P_L(r_1,r_2,r_3)| \sim \frac{8(e^{\epsilon/2}-e^{-\epsilon/2})^3  }{ vol(T_1(S))|Isom^+(r_1,r_2,r_3)|}e^{(r_1+r_2+r_3)L/2}
\end{eqnarray*}
where $Isom^+(r_1,r_2,r_3)$ is the orientation-preserving isometry group of the pair of pants with boundary lengths $r_1,r_2,r_3$.
\end{cor}

For comparison, recall that the number of closed oriented geodesics in $S$ with length in $(L-\epsilon,L+\epsilon)$ is asymptotic to $(e^{\epsilon}-e^{-\epsilon})e^L/L$ (see e.g. \cite{Bus}). Unlike $|\P_L(r_1,r_2,r_3)|$ it does not depend on the genus of the surface. 

We will prove this as a corollary to theorem \ref{thm:counting} below. Recall that a closed oriented geodesic $\gamma$, is a local isometry $\gamma: S^1 \to S$ from the circle of length $length(\gamma)$ to $S$. We identify geodesics $\gamma_1,\gamma_2$ if there is an orientation-preserving isometry $\Psi: S^1 \to S^1$ such that $\gamma_1 = \gamma_2 \circ \Psi$.  Let $\bgamma$ denote the image of $\gamma$ so $length(\bgamma) = length(\gamma)/m$ where the map $\gamma: S^1 \to \bgamma$ is an $m$-fold cover.

For $\gamma \in \G_{r_1L}=\G_{r_1L}(-\epsilon,\epsilon)$, let $\P_L(r_1,r_2,r_3;\gamma) \subset \P_L(r_1,r_2,r_3)$ denote the subset of immersions $(j: P \to S)$ in which $j$ restricted to some boundary component is equivalent to $\gamma$.

\begin{thm}\label{thm:counting}
Assume $Comm(\Gamma)<Isom^+(\H^2)$. Let $\gamma_L \in \G_{r_1L}$. If $r_2+r_3 > r_1$ and $r_1 + r_3 > r_2$ then 
\begin{eqnarray*}
|\P_L(r_1,r_2,r_3;\gamma_L)| \sim \frac{4(e^{\epsilon/2}-e^{-\epsilon/2})^2 }{ vol(T_1(S)) n} length(\bgamma_L) \exp(-length(\gamma_L)/2 +r_2L/2 + r_3L/2)
\end{eqnarray*}
where $n=1$ if $r_2 \ne r_3$, otherwise $n=2$.
\end{thm}
Note if $\P_L(1,1,1;\gamma)=\P_L(1,1,1;-\gamma)$ for all $\gamma \in \G_L$ (where $-\gamma$ equals $\gamma$ with reversed orientation) then question \ref{cnj:singlesurface} can be answered affirmatively by gluing the immersions in $\P_L$ together in such a way that every pair of pants $P \in \P_L$ is used exactly once.

% An integral formula for $f(\epsilon)$ is given at the end of the proof. It seems possible to obtain a closed form for $f$ from it.

%----------new stuff--------------

By successively gluing together pants immersions, the above theorem can be used to give asymptotic counts for immersions of $n$-holed spheres. For example, suppose $T$ is an $n$-holed sphere with a fixed pants decomposition. Consider locally isometric immersions $j: {\tilde T} \to S$ where ${\tilde T}$ equals $T$ with a hyperbolic metric in which the boundary is totally geodesic and all curves of the pants decomposition have length in $L+(-\epsilon,\epsilon)$. Then the number of such immersions is asymptotic to
\begin{eqnarray*}
\frac{h(\epsilon)L^{n-3}e^{nL/2}}{vol(T_1(S))^{n-2}}
\end{eqnarray*}
where $h(\epsilon)\approx \epsilon^{2n-3}$ depends on the symmetries of the decomposition.

%By contrast, the dimension of the moduli space of $n$-holed spheres is $3n-6$.

The theorem above can be improved by taking ``twisting'' into
account. For $I \subset \bgamma$ let $\P_L(r_1,r_2,r_3;\gamma,I)
\subset \P_L(r_1,r_2,r_3;\gamma)$ be the subset of pants immersions $j: P
\to S$ such that the shortest curve $\beta \subset P$ from $\partial_1 P$ to
$\partial_2 P$ has an endpoint on $\partial_1 P$ which maps to
$I$. Here $j$ restricted to $\partial_1 P$ equals $\gamma$ and $j$
restricted to $\partial_2 P$ is a closed geodesic of length in $r_2L +
(-\epsilon,\epsilon)$.

Define $\tM=\tM(L)$ by:
\begin{displaymath}
\tM(L)=\left\{\begin{array}{ll}
    e^{(-r_1 -r_2+ r_3)L/4} & \textnormal{ if $r_1,r_2,r_3$ are the
    sidelengths of a nondegenerate triangle}\\
    1 & \textnormal{ if $r_1+r_2=r_3$}\\
    (-r_1-r_2 + r_3)L/2  & \textnormal{ if $r_3 > r_1 + r_2$}
\end{array}\right.
\end{displaymath}

\begin{thm}\label{thm:twist}
Assume $Comm(\Gamma)<Isom^+(\H^2)$, $r_2+r_3 > r_1$ and $r_1 + r_3 >
r_2$. Let $\gamma_L \in \G_{r_1L}$ and $I_L$ be a subsegment of
$\bgamma_L$. If $length(I_L) \sinh(\tM)e^{r_2L/2}$ is bounded away from $0$ then
\begin{eqnarray*}
|\P_L(r_1,r_2,r_3;\gamma_L,I_L)| \sim \frac{length(I_L) |\P_L(r_1,r_2,r_3;\gamma_L)|}{length(\bgamma)}.
\end{eqnarray*}
\end{thm}
The hypothesis on $length(I_L)$ is automatically satisfied if
$length(I_L)$ is constant. The coefficient $\sinh(\tM)e^{r_2L/2}$ is
probably not optimal. $\tM(L)$ is roughly the length of the shortest path from $\partial_1 P$ to $\partial_2 P$. From this result, the Ehrenpreis conjecture follows if we weaken the conclusion so that ${\tilde S_i}$ is no longer required to be closed but is required to have maximum injectivity radius bounded from above (see corollary \ref{cor:ehrenpreis}).

An application: let $T$ be a hyperbolic 4-holed sphere with its boundary components labeled $c_1,c_2,c_3,c_4$. Let $c_5$ be an oriented simple closed geodesic in the interior separating $c_1$ and $c_2$ from $c_3$ and $c_4$. Let $\alpha$ be the shortest path from $c_1$ to $c_5$ and let $\beta$ be the shortest path from $c_5$ to $c_3$. Let $\tau$ be the twist parameter of $c_5$: $\tau$ equals the signed distance from the endpoint of $\alpha$ to the endpoint of $\beta$ along $c_5$ with sign determined the orientation of $c_5$. The hyperbolic structure on $T$ is determined by the lengths of $c_1,..,c_5$ and $\tau$. Given $r_1,r_2,r_3,r_4, r_5>0$ the number of immersions $T \to S$ in which $length(c_i) \in r_i L + (-\epsilon,\epsilon)$ for all $i$ and $|\tau| < \epsilon$ is asymptotic to:
\begin{eqnarray*}
\frac{h(\epsilon)e^{(r_1+r_2+r_3+r_4)L/2}}{vol(T_1(S))^2}
\end{eqnarray*}
where $h(\epsilon)\approx \epsilon^6$. The formula does not depend on $r_5$. One explanation for this is that there are three simple closed curves in the interior of a 4-holed sphere. If, for example, we required $c_5$ to separate $c_1,c_3$ from $c_2,c_4$ instead then its length might change but the asymptotic count cannot. 

The techniques of this paper can also be used to obtain asymptotic counts
for immersions of bordered surfaces with positive genus although that
is not taken up here.

%------------new stuff----------------

\subsection{Clotheslines}

The next results have to do with the ``distance'' between geodesics in $\G_L=\G_L(-\epsilon,\epsilon)$. By abuse of notation, we identify the pants $P$ with the immersion $j:P\to S$. For example, we may write $P \in \P_L(r_1,r_2,r_3)$ to mean $(j:P \to S) \in \P_L(r_1,r_2,r_3)$.

\begin{defn}
For $\alpha, \beta \in \G_L$, an $L$-clothesline from $\alpha$ to $\beta$ is an $(2n-1)$-tuple $(\alpha=\gamma_1,P_1,\gamma_2,P_2,...,\gamma_{n-1},P_n,\gamma_n=\beta)$ where $\gamma_i \in \G_L$ for all $i$, and $P_i \in \P_L(1,1,1)$ has boundary components mapping to $\gamma_{i-1}$ and $-\gamma_i$ where $-\gamma_i$ denotes $\gamma_i$ with orientation reversed. It is possible to glue the immersions $j_i:P_i \to S$ together to obtain an immersion
\begin{eqnarray*}
j: P_1 \cup_{\gamma_2} P_2 \cup_{\gamma_3} ... \cup_{\gamma_{n-1}} P_n \to S
\end{eqnarray*}
from an $(2+n)$-holed sphere to $S$. The length of the clothesline is $n$. If $\gamma \in \G_L$ has the property that for every $\gamma' \in \G_L$ there is an $L$-clothesline of length $2$ from $\gamma$ to $\gamma'$ then we say that $\gamma$ is regular.
\end{defn}

\begin{thm}\label{thm:clotheslines}
Assume $Comm(\Gamma)<Isom^+(\H^2)$. Then for all $L$ sufficiently large, for any pair $\gamma_1,\gamma_2 \in \G_L$ there is an $L$-clothesline of length 3 from $\gamma_1$ to $\gamma_2$. Also
\begin{eqnarray*}
\lim_{L\to\infty} \frac{\#\{\textnormal{regular geodesics in }\G_L\}}{\#\G_L} = 1.
\end{eqnarray*} 
\end{thm}

If every geodesic in $\G_L$ is regular, then question \ref{cnj:singlesurface} can be answered affirmatively. To see this, let $\gamma_1, \sigma_1 \in \G_L$. Assuming $\gamma_1$ say is regular, there exists an
immersed four-holed sphere $H_1$ in $S$ obtained from
gluing two pants immersions $P_1, P_2 \in \P_L(1,1,1)$ together so that two of the boundary
components of $H_1$ are $\gamma_1$ and $\sigma_1$. Let $\gamma_2$ and
$\sigma_2$ be the other two boundary components of $H_1$. If $\gamma_2$ say is regular, there exists an immersed four-holed sphere $H_2$ which decomposes as the
union of two pants in $\P_L(1,1,1)$ and such that $\gamma_2$ and
$\sigma_2$ are two of its boundary components. Suppose this process can be continued indefinitely. Since $\G_L$ is finite, at some later stage we must use the same pair of geodesics in $\G_L$ that we used at an earlier
stage. That is,
$(\gamma_m, \sigma_m)= (\gamma_n,\sigma_n)$ for some $m < n$. Let
$H_m, H_{m+1},.., H_{n-1}$ be the intervening four-holed spheres. By
construction, it is possible to glue $H_m$ to $H_{m+1}$ to ... to
$H_{n-1}$ to $H_m$ to obtain a closed surface ${\tilde S}$ with a
locally isometric immersion into $S$ answering question
\ref{cnj:singlesurface}.

\subsubsection{Organization of Paper}

\begin{itemize}

\item section \ref{section:equidistribution}:  The equidistribution results described in the introduction are proven here.  

\item section \ref{section:isometry}: The main technique for constructing elements $g \in \Gamma$ by perturbations is detailed here.

\item section \ref{section:pants}: Methods of the previous section are used to construct pants immersions.

\item section \ref{section:counting}: The asymptotic counting results presented in the introduction are proven here. 

\item section \ref{section:clotheslines}: The clothesline results are proven here.

\item section \ref{section:calculations}: Main calculations used throughout the paper are presented here. In particular, corollary \ref{cor:perp} is proven.

\end{itemize}

{\bf Acknowledgements}: I'd like to thank Joel Hass with whom conversations about the surface subgroup conjecture eventually led to this project. I'd also like thank Chris Connell and Chris Judge for many helpful conversations.

\section{Notation}\label{section:notation}

%This section has been checked 4/16/05. pants16.tex.

Throughout the paper, $\Gamma <PSL_2(\R)$ will denote a fixed cocompact discrete group. Let $S=\H^2/\Gamma$. Let $w_0=(i,i)$ be the unit vector in the upperhalf plane model of $\H^2$. Then $PSL_2(\R)$ is identified with the unit tangent bundle $T_1(\H^2)$ by $g \to gw_0$. This identification projects to an identification between $T_1(S)$ and $\Gamma \backslash PSL_2(\R)$. If $w \in T_1(\H^2)$ and $w=g w_0$ for $g \in PSL_2(\R)$ we may write $w^{-1}$ to mean $g^{-1}$ or $wh$ to mean $gh$ when $h \in PSL_2(\R)$. 

A pair of pants $P$ is a compact hyperbolic surface with geodesic boundary homeomorphic to the $2$-sphere minus 3 open disks. An immersion $j: P \to S$ will always mean a locally-isometric immersion. Define:

\begin{displaymath}\begin{array}{l}
V_T:=\left[ \begin{array}{cc}
1 & 0 \\
T & 1
\end{array} \right] \,
U_T:=\left[ \begin{array}{cc}
1 & T \\
0 & 1
\end{array} \right]\,
R_\delta:=
\left[ \begin{array}{cc}
\cos(\delta/2) & \sin(\delta/2) \\
-\sin(\delta/2) & \cos(\delta/2)
\end{array} \right] \\
\\
G_L:=\left[ \begin{array}{cc}
e^{L/2} & 0 \\
0 & e^{-L/2} 
\end{array} \right] \,
Mat[a,b,c,d]:=
\left[ \begin{array}{cc}
a & b \\
c & d
\end{array} \right].
\end{array}
\end{displaymath}

We identify these matrices with their images in $PSL_2(\C)$. The action $\Gamma g \to \Gamma g G_t$ on $\Gamma \backslash PSL_2(\R)$ is the geodesic flow for time $t$. We identify these matrices with their images in $PSL_2(\C)$. The actions $\Gamma g \to \Gamma g U_t$ and $\Gamma g \to \Gamma g V_t$ are distinct horocycle flows (for time $t$). $R_\delta$ is a rotation by the angle $\delta$ about the point $i$ in the upperhalf plane model. When discussing a matrix of the form $Mat[a,b,c,d]$ unless otherwise specified we assume $a>0$ and $ad-bc=1$.

%Let $\sB_\epsilon =\{Mat[a,b,c,d] \in PSL_2(\R): |\ln(a)|,|b|,|c|,|d-1| < \epsilon, a>0 \}$. $\sB_\epsilon$ is a neighborhood of the identity.

\section{Equidistribution Theorems}\label{section:equidistribution}

In this section, theorem \ref{thm:equidistribution}, theorem \ref{thm:angle} and a special case of Rufus Bowen's equidistribution theorem (theorem \ref{thm:bowen}) are proven. But first, we note that the hypothesis $Comm(\Gamma)< Isom^+(\H^2)$ is generic.

\begin{remark}\label{rmk:margulis}
\textnormal{ By work of Margulis (see e.g. \cite{Zim}), if $\Gamma$ is not arithmetic, then $\Gamma$ has finite index in $Comm(\Gamma)$. Hence $\H^2/\Gamma$ finitely covers $\H^2/Comm(\Gamma)$. The covering map induces a map from the moduli space of $\H^2/Comm(\Gamma)$ into the moduli space of $\H^2/\Gamma$. It can be shown that if $Comm(\Gamma)$ contained orientation-reversing elements then the dimension of the moduli space of $\H^2/Comm(\Gamma)$ is strictly less that the dimension of the moduli space of $\H^2/\Gamma$. Thus, $Comm(\Gamma)<Isom^+(\H^2)$ holds generically.}\end{remark}

\subsection{Unipotent Flows on the Product}
Consider a unipotent flow on supported on 
\begin{eqnarray*}
\Gamma \backslash PSL_2(\R) \times \Gamma \backslash PSL_2(\R) = T_1(S) \times T_1(S)
\end{eqnarray*}
of the form $(\Gamma g_1,\Gamma g_2) \to (\Gamma g_1 U_t, \Gamma g_2 kU_tk^{-1})$ where $k \in Isom(\H^2)$ and $U_t$ is given in section \ref{section:notation}. For $g_1,g_2 \in PSL_2(\R)$ let $\omega_{T,g_1,g_2}$ be the probability measure on 
\begin{eqnarray*}
\{(\Gamma g_1 U_t,\Gamma g_2 kU_tk^{-1}): t\in [0,T]\}
\end{eqnarray*}
obtained from Lebesgue measure on $[0,T]$. The orbit of $(\Gamma g_1,\Gamma g_2)$ under the flow is said to be equidistributed with respect to a probability measure $\lambda'$ on $T_1(S) \times T_1(S)$ if $\omega_{T,g_1,g_2}$ converges to $\lambda'$ in the weak* topology as $T$ tends to infinity. The orbit is said to be exceptional if it not equidistributed with respect to the product Haar measure $\lambda \times \lambda$. The next lemma implies that the commensurator of $\Gamma$ essentially classifies the exceptional orbits.

\begin{lem}\label{lem:Ratner}
Suppose that the orbit of $(\Gamma g_1,\Gamma g_2)\in T_1(S) \times T_1(S)$ under the unipotent flow $(\Gamma g_1,\Gamma g_2) \to (\Gamma g_1 U_t, \Gamma g_2 kU_tk^{-1})$ is not equidistributed with respect to the product Haar measure on $T_1(S) \times T_1(S)$. Then there exists a closed subgroup $H < PSL_2(\R) \times PSL_2(\R)$ such that
\begin{itemize}
\item the closure of the orbit of $(\Gamma g_1,\Gamma g_2)$ equals $(\Gamma g_1, \Gamma g_2)H$,

\item $Stab(\Gamma g_1,\Gamma g_2) \cap H$ is a lattice in $H$ where $Stab(\Gamma g_1,\Gamma g_2)$ is the stabilizer of $(\Gamma g_1,\Gamma g_2)$ in $PSL_2(\R) \times PSL_2(\R)$,

\item the orbit of $(\Gamma g_1,\Gamma g_2)$ is equidistributed with respect to the probability measure $\lambda_H$ on $(\Gamma g_1,\Gamma g_2)H$ induced by Haar measure on $H$,

\item $H=H_T:=\{(g, k U_T g U_T^{-1} k^{-1}): g \in PSL_2(\R)\}$ for some $T \in \R$,

\item $g_2 k U_{T} g_1^{-1} \in Comm(\Gamma)$.
\end{itemize}
\end{lem}

%need to explain the abuse of notation from T_1(S) to \Gamma \backslash PSL_2(\R).

\begin{proof}
The first three statements follow from Ratner's work on Ragnuthan's conjectures \cite{Ra5}. Both factors of $(\Gamma \times \Gamma)H$ project onto $\Gamma \backslash PSL_2(\R)$. The classification of closed subgroups of $PSL_2(\R) \times PSL_2(\R)$ implies that $H$ has the form
\begin{eqnarray*}
H= \{(g, hgh^{-1} | g\in PSL_2(\R)\}
\end{eqnarray*}
for some $h \in Isom(\H^2)$. The fourth statement now follows from the fact that the group $\{(U_t,k U_t k^{-1}): t \in \R\}$ fixes $(\Gamma \times \Gamma)H$. 

We show that the last statement follows from the second and fourth. Note
\begin{eqnarray*}
(g_1,g_2)^{-1}(\Gamma \times \Gamma )(g_1,g_2)=Stab(\Gamma g_1,\Gamma g_2).
\end{eqnarray*}
Hence $(\Gamma \times \Gamma) \cap (g_1,g_2) H (g_1,g_2)^{-1}$ is a lattice in $(g_1,g_2) H (g_1,g_2)^{-1}$. But
\begin{eqnarray*}
(g_1,g_2) H (g_1,g_2)^{-1} &=& \Big\{\big(g_1gg_1^{-1},g_2 k U_T gU_T^{-1}k^{-1} g_2^{-1}\big) : g\in PSL_2(\R)\Big\}\\
                   &=& \Big\{ \big(g, (g_2k U_T g_1^{-1}) g (g_1 U_T^{-1} k^{-1} g_2^{-1})\big) : g \in PSL_2(\R)\Big\}.
\end{eqnarray*}
Thus $g_2k U_T g_1^{-1} \in Comm(\Gamma)$.
\end{proof}

\begin{cor}
If $Comm(\Gamma)$ does not contain any orientation-reversing elements, then the orbit of any $(\Gamma g_1,\Gamma g_2) \in \Gamma \backslash PSL_2(\R) \times \Gamma \backslash PSL_2(\R)$ under the flow $(\Gamma  g_1,\Gamma g_2) \to (\Gamma g_1U_t,\Gamma g_2V_t)$ is equidistributed with respect to the product Haar measure on $T_1(S)\times T_1(S)$ where $V_t$ is defined in section \ref{section:notation}.
\end{cor}

\begin{proof}
There exists $k \in Isom(\H^2)$ with $V_t = k U_t k^{-1}$. Every such $k$ is orientation-reversing. Hence $g_2kU_Tg_1^{-1}$ is orientation-reversing for any $T\in\R, g_1,g_2 \in \Gamma$. The fifth statement of the lemma above now implies the corollary.
\end{proof}

\subsection{Hypercycle Flows}
\begin{figure}[htb] 
\begin{center}
 \ \psfig{file=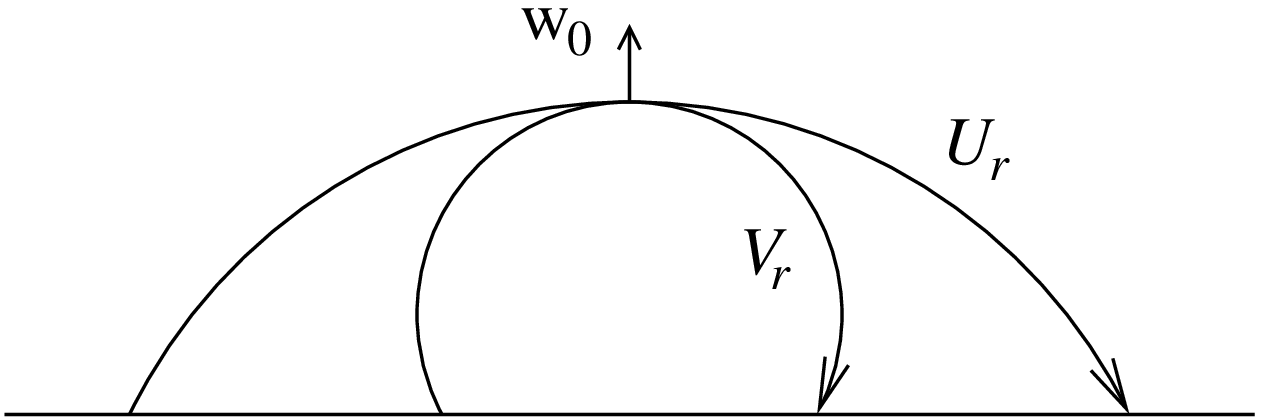,height=1in,width=3in}
 \caption{Hypercycles $\sV_r$ and $\sU_r$.}
 \label{fig:hypercycle}
 \end{center}
 \end{figure}

Given a geodesic $\sigma \subset \H^2$, its radius $r$-neighborhood is bounded by two constant-curvature curves which we call hypercycles of radius $r$. The unit vector $w_0=(i,i)$ (in the upperhalf plane model) is normal to two distinct radius $r$ hypercycles, $\sV_r$ and $\sU_r$. Assume that $\sV_r$ partially bounds the radius $r$ neighborhood $N_r(\sigma)$ of a geodesic for which $w_0$ is an outer unit normal. See figure \ref{fig:hypercycle}. Orient $\sV_r$ and $\sU_r$ so that their endpoints at infinity in the upper half plane model are ordered from negative to positive. Define $w_0 V_{r,t} \in T_1(\H^2)$ so that $w_0 V_{r,t}$ is normal to $\sV_r$, on the same side of $\sV_r$ as $w_0$ and the map $t \to basepoint(w_0 V_{r,t})$ parametrizes $\sV_r$ by arclength in an orientation-preserving manner. This defines a flow on $T_1(\H^2)$ by setting $(gw_0) V_{r,t}= g(w_0V_{r,t})$ for any $g \in PSL_2(\R)$. Define $U_{r,t}$ similarly. $V_{r,t}$ and $U_{r,t}$ are distinct hypercycle flows that converge to $V_t$ and $U_t$ as $r \to \infty$. These flows descend to flows on $\Gamma\backslash PSL_2(\R)$.

\begin{proof} (of theorem \ref{thm:equidistribution})
For $T>0$ and $z=(z_1,z_2) \in T_1(S)$, let $\nu_{T,z}$  be the probability measure on 
\begin{eqnarray*}
\{(z_1 U_t,z_2 V_t): t\in [0,T]\}
\end{eqnarray*}
obtained from Lebesgue measure on $[0,T]$. By the previous corollary $\nu_{T,z}$ converges to $\lambda \times \lambda$ in the weak* topology as $T \to \infty$. Convergence is uniform in $z$ since $T_1(S) \times T_1(S)$ is compact. Thus there exists a $T'_0>0$ such that for all $T> T'_0$, 
\begin{eqnarray*}
|\nu_{T,z}(f)-\lambda \times\lambda(f)|<\epsilon_0/2
\end{eqnarray*}
for all $z$. A standard calculation shows that for any $w\in T_1(S)$ if $w_t:=w R_{-\pi/2} G_t R_{\pi/2}$ then
\begin{eqnarray*}
(w_t G_{-L/2}, w_t G_{L/2}) = \Big((w G_{-L/2})U_{L/2,t\cosh(L/2)}, (w G_{L/2})V_{L/2,t\cosh(L/2)}\Big).
\end{eqnarray*}
%we have checked the above.
Thus $\mu_{w,T,L}$ is the probability measure on 
\begin{eqnarray*}
\Big\{\big( (w G_{-L/2})U_{L/2,t\cosh(L/2)e^{-L/2}}, (w G_{L/2})V_{L/2,t\cosh(L/2)e^{-L/2}}\big) \, \big| 0\le t \le T\Big\}
\end{eqnarray*}
induced from Lebesgue measure on $[0,T]$. Since $(U_{L/2,t}, V_{L/2,t})$ converges to $(U_t,V_t)$ as $L \to \infty$, it follows that any weak* subsequential limit of $\{\mu_{w,T,L}\}_{L >0}$ is equal to $\nu_{T/2,z}$ for some $z \in T_1(S)\times T_1(S)$. Thus there exists a $T_0 \ge T_0'$ such that
\begin{eqnarray*}
|\mu_{w,T,L}(F)- \lambda\times\lambda(F)| < \epsilon_0
\end{eqnarray*}
 for all $L> T_0$.
\end{proof}

%\begin{remark}
%For a generic uniform lattice $\Gamma < PSL_2(\R)$, $Comm(\Gamma)$ has no orientation-reversing elements.
%\end{remark}

\subsection{Proof of theorem \ref{thm:angle}}

To prove theorem \ref{thm:angle} we need the following calculation which is proven as a corollary to theorem \ref{thm:main} below.
\begin{cor}\label{cor:perp}
Let $F=(F_1,F_2,F_3):(0,\infty) \times \R^2 \to \R^3$ be defined by
\begin{eqnarray*}
F(a,b,c)=\Big(\frac{c-b}{a}, \frac{c+b}{a}, 2\ln(a) \Big).
\end{eqnarray*}
Let $B>0$. Suppose $w \in T_1(\H^2)$, $g\in PSL_2(\R)$ and $gwG_{-L/2}=wG_{L/2}Mat[a,b,c,d]$ for some $a,b,c,d \in \R$ with $|\ln(a)|,|b|,|c|,|d-1|<B$ and $ad-bc=1$. Let $axis(w)$ denote the geodesic in $\H^2$ tangent to $w$ oriented consistently with $w$. Let $w^\perp$ be a unit vector perpendicular to $w$ so that $(w,w^\perp)$ is a positively oriented frame. Then the signed distance between the basepoint of $w$ and $axis(w^\perp) \cap axis(g)$ equals 
\begin{eqnarray*}
F_1(a,b,c)e^{-L/2} + O((B+1)^3 e^{-L}).
\end{eqnarray*}
The sign is positive iff $axis(w^\perp)\cap axis(g)$ is on the right side of $axis(w)$. The angle from $axis(g)$ to $axis(w^\perp)$ equals 
\begin{eqnarray*}
\pi/2 + F_2(a,b,c)e^{-L/2} + O((B+1)^2 e^{-L}).
\end{eqnarray*}
 Also the translation length $g$ satisfies
\begin{eqnarray*}
tr.length(g)= L + F_3(a,b,c) + O((B+1)e^{-L}).
\end{eqnarray*}
\end{cor}

We will also need the following volume identity.

\begin{lem}\label{lem:vol}
For a compact set $X \subset (0,\infty) \times \R^2$ define
\begin{eqnarray*}
\B_X = \{Mat[a,b,c,d]\in PSL_2(\R): F(a,b,c)\in X\}
\end{eqnarray*}
where $F(\cdot)$ is defined in corollary \ref{cor:perp}. Then the volume, $vol(\B_X)$, of $\B_X$ with respect to the standard volume form on $PSL_2(\R)$ satisfies
\begin{eqnarray*}
vol(\B_X) = \int_X e^{x_3}  dx
\end{eqnarray*}
where $x=(x_1,x_2,x_3)$ and the integral is with respect to Euclidean volume.
\end{lem}

\begin{proof}
Without loss of generality, we may assume that $X$ is the product of compact intervals $X=I_1 \times I_2 \times [l_1,l_2]$ and that $X$ is small enough so that for any $g\in PSL_2(\R)$ the map $h \to \Gamma gh \in \Gamma \backslash PSL_2(\R)$ is injective when restricted to $\B_X$. 

Suppose $\tv$ is a unit vector in $T_1(\H_2)$ and there exists an element $g \in \Gamma$ such that
\begin{eqnarray*}
g \tv G_{-L/2}\in \tv G_{L/2} \B_X.
\end{eqnarray*}
Then $g$ is unique by hypotheses on $X$, so we may write $g=g_{\tv}$. The projection of $axis(g_{\tv})$ to $S$ is an oriented geodesic $\gamma \in \G_L(l_1,l_2)$. Let $v$ be the projection of $\tv$ to $T_1(S)$. $\gamma$ is uniquely determined by $v$, so we may write $\gamma=\gamma_v$. 

Let $\Omega_L \subset T_1(S)$ be the set of all vectors $v$ such that $\gamma_v$ exists and is in $\G_L(l_1,l_2)$. The proof proceeds by calculating the volume of $\Omega_L$ in two different ways.

For any geodesic $\gamma \in \G_L(l_1,l_2)$ define
\begin{eqnarray*}
\Omega(\gamma)=\{v \in T_1(S): \gamma = \gamma_v\}.
\end{eqnarray*}
From the corollary \ref{cor:perp}, it follows that $\Omega(\gamma)$ is approximately the set of all vectors $v \in T_1(S)$ such that there is a geodesic segment $\beta \subset S$ satisfying 
\begin{itemize}
\item $v$ is based at an endpoint of $\beta$ and is perpendicular to $\beta$ there,
\item the other endpoint of $\beta$ is in $\gamma$
\item $length(\beta) \in I_1e^{-L/2} + O(e^{-L})$
\item the angle from $\beta$ to $\gamma$ is in $I_2e^{-L/2} + O(e^{-L})$.
\end{itemize}
Hence the volume of $\Omega(\gamma)$ satisfies
\begin{eqnarray*}
vol(\Omega(\gamma)) = length(\gamma)length(I_1)length(I_2) e^{-L} + O(Le^{-3L/2}).
\end{eqnarray*}
The hypotheses on $X$ imply that for different $\gamma_1,\gamma_2 \in \G_L(l_1,l_2)$, the intersection $\Omega(\gamma_1)\cap \Omega(\gamma_2)$ is empty. By definition,
\begin{eqnarray*}
\Omega_L = \bigcup_{\gamma \in \G_L(l_1,l_2)} \Omega(\gamma).
\end{eqnarray*}
So
\begin{eqnarray*}
vol(\Omega_L) &=& \Sigma_{\gamma \in \G_L(l_1,l_2)} vol(\Omega(\gamma))\\
              &=&  \Sigma_{\gamma \in \G_L(l_1,l_2)} \big[length(\gamma)length(I_1)length(I_2) e^{-L} + O(Le^{-3L/2})\big].
\end{eqnarray*}
It is well-known that 
\begin{eqnarray*}
\# \G_L(l_1,l_2) \sim \frac{e^{l_2}-e^{l_1}}{L}e^{L}.
\end{eqnarray*}
See for example \cite{Bus}. Thus 
\begin{eqnarray*}
vol(\Omega_L) &\sim& length(I_1)length(I_2)(e^{l_2}-e^{l_1}).
\end{eqnarray*}
Let $f_L:T_1(S) \to T_1(S)\times T_1(S)$ be the map
\begin{eqnarray*}
f_L(v)=(v G_{-L/2}, v G_{L/2}).
\end{eqnarray*}
Let $\lambda$ be Haar probability measure on $T_1(S)$. The key observation is that the equidistribution theorem \ref{thm:equidistribution} implies that the push-forward measure $f_{L*}(\lambda)$ converges to $\lambda \times \lambda$ in the weak* topology as $L$ tends to infinity. So
\begin{eqnarray*}
vol(\Omega_L) &=& vol(T_1(S)) \lambda(\Omega_L)\\
                 &=& vol(T_1(S)) f_{L*}(\lambda)(\{(\Gamma g_1, \Gamma g_2) \, | \, \Gamma g_1 \in \Gamma g_2 \sB_X \}) \\
                 &\sim& vol(T_1(S)) \lambda \times \lambda (\{(\Gamma g_1, \Gamma g_2) \, | \, \Gamma g_1 \in \Gamma g_2 \sB_X \}) \\
                 &=& vol(\B_X).
\end{eqnarray*}
Hence $vol(\B_X)= length(I_1)length(I_2)(e^{l_2}-e^{l_1}) = \int_X e^{x_3} dx$.

\end{proof}

Now we are ready to prove theorem \ref{thm:angle}
\begin{proof}(of theorem \ref{thm:angle})
Let $I$ be a finite interval in $\R$. Let $X=I \times [a_1,a_2] \times [l_1,l_2]$, let $\B_X$ be as defined in the previous lemma. Without loss of generality we may assume that  $X$ is small enough so that for any $g\in PSL_2(\R)$ the map $h \to \Gamma gh \in \Gamma \backslash PSL_2(\R)$ is injective when restricted to $\B_X$. 

Let $W_L$ be the set of all unit vectors $v$ perpendicular to $\sigma_L$. To prove the theorem, we will estimate the length of $W_L \cap \Omega_L$ in two different ways where $\Omega_L$ is as defined in the previous lemma.

It follows from corollary \ref{cor:perp} that each connected component of $\Omega_L \cap W_L$ has length equal to $length(I_1)e^{-L/2} + O(e^{-L})$ except if this component contains an endpoint of $\sigma_L$. The latter contributes an asymptotically negligible amount to the length of $W_L \cap \Omega_L$ and so it can be ignored. Thus,
\begin{eqnarray*}
length(W_L \cap \Omega_L) \sim \#N_L(a_1,a_2,l_1,l_2)length(I_1)e^{-L/2}.
\end{eqnarray*}
Let $w \in W_L$ be based at an endpoint of $\sigma_L$. Let $\mu = \mu_{\tw,T,L}$ be as defined in theorem \ref{thm:equidistribution} where $\tw \in T_1(\H^2)$ projects to $w$. Let $T_L = e^{L/2}length(\sigma_L)$. It follows from that theorem that
\begin{eqnarray*}
\frac{length(W_L \cap \Omega_L)}{length(\sigma_L)} &=& \pi_*(\mu_{\tw, T_L,L})(\{(\Gamma g_1, \Gamma g_2) \, | \, \Gamma g_1 \in \Gamma g_2 \sB_X \}) \\
&\sim& \lambda \times \lambda (\{(\Gamma g_1, \Gamma g_2) \, | \, \Gamma g_1 \in \Gamma g_2 \sB_X \}) \\
&=& \frac{vol(\sB_X)}{vol(T_1(S))}.
\end{eqnarray*}
So
\begin{eqnarray*}
\#N_L(a_1,a_2,l_1,l_2)length(I_1)e^{-L/2}&\sim&\frac{vol(\sB_X)length(\sigma_L)}{vol(T_1(S))}.
\end{eqnarray*}
The previous lemma now implies the theorem.

\end{proof}

\subsection{Rufus Bowen's Equidistribution Theorem}

For any closed geodesic $\gamma \subset T_1(S)$ let $\nu_\gamma$ be the measure on $T_1(S)$, supported on the vectors tangent to $\gamma$, induced by Lebesgue measure. The total mass of $\nu_\gamma$ equals the length of $\gamma$. Recall that $\G_L$ is the set of all oriented closed geodesics $\gamma$ in $S$ with length in $(L-\epsilon,L+\epsilon)$. 

The next theorem was proven first by Rufus Bowen \cite{RBow} in the more general context of Axiom A flows. Our proof is significantly shorter although we invoke Ratner's theorems (lemma \ref{lem:Ratner}). 
\begin{thm}\label{thm:bowen}
Assume $Comm(\Gamma)<Isom^+(\H^2)$. Given $L,\epsilon$, let
\begin{eqnarray*}
\nu_{L} = \nu_{L,\epsilon} = \frac{1}{N}\Sigma_{\gamma \in \G_L} \, \nu_\gamma
\end{eqnarray*}
where $N>0$ is chosen so that $\nu_{L,\epsilon}$ is a probability measure. Then, as $L$ tends to infinity, $\nu_{L,\epsilon}$ converges to $\lambda$, the Haar probability measure on $T_1(S)$, in the weak* topology.
\end{thm}

\begin{proof}
Let $X=[-\epsilon,\epsilon]^3$. Let $\B_X$, $\Omega(\gamma),\Omega_L$ be as defined in lemma \ref{lem:vol}. Let $\omega_\gamma$ be the measure on $T_1(S)$ defined by
\begin{eqnarray*}
\omega_\gamma = \chi_{\Omega(\gamma)} \lambda
\end{eqnarray*}
where $\chi_{\Omega(\gamma)}$ is the characteristic function of
$\Omega(\gamma)$. The total mass of $\omega_\gamma$ is the
$\lambda$-measure of $\Omega(\gamma)$. Let
\begin{eqnarray*}
\omega_{L} = \frac{1}{\lambda(\Omega_L)}\Sigma_{\gamma \in \G_L} \, \omega_\gamma = \frac{1}{\lambda(\Omega_L)} \chi_{\Omega_L}\lambda.
\end{eqnarray*}
$\omega_L$ is a probability measure.

Claim: $\omega_L$ converges to $\lambda$ in the weak* topology as $L$ tends to infinity.

Proof:  Let $f_L:T_1(S) \to T_1(S)\times T_1(S)$ be the map
\begin{eqnarray*}
f_L(w)=(w G_{-L/2}, w G_{L/2}).
\end{eqnarray*}
Let $W$ be an open subset of $T_1(S)$. Let $\lambda_W$ be the probability measure 
\begin{eqnarray*}
\lambda_W = \frac{\chi_W \lambda}{\lambda(W)}.
\end{eqnarray*}
The key observation is that the equidistribution theorem \ref{thm:equidistribution} implies that the push-forward measure $f_{L*}(\lambda_W)$ converges to $\lambda \times \lambda$ in the weak* topology as $L$ tends to infinity. So
\begin{eqnarray*}
\omega_L(W) &=& \frac{\lambda(\Omega_L \cap W)}{\lambda(\Omega_L)}\\
                 &=& \frac{\lambda_W(\Omega_L)\lambda(W)}{\lambda(\Omega_L)}\\
                 &=& \frac{f_{L*}(\lambda_W)(\{(\Gamma g_1, \Gamma g_2) \, | \, \Gamma g_1 \in \Gamma g_2 \sB_X\})\lambda(W)}{\lambda(\Omega_L)}\\
                 &\sim& \frac{\lambda(\Gamma \sB_X)\lambda(W)}{\lambda(\Omega_L)}.
\end{eqnarray*}
Substitute $W=T_1(S)$ in the above to obtain 
\begin{eqnarray*}
1 \sim \frac{\lambda(\Gamma \sB_X)}{\lambda(\Omega_L)}.
\end{eqnarray*}
Thus $\omega_L(W) \sim \lambda(W)$. Since $\Omega$ is arbitrary, $\omega_L$ converges to $\lambda$ as claimed. 

Now let $f:T_1(S)\to \C$ be any continuous function. We will show that $|\omega_L(f)-\nu_L(f)|$ tends to zero as $L$ tends to infinity. Since $\omega_L$ converges to $\lambda$ this will imply the theorem.

Let $A_L=\epsilon^2 e^{-L}$. It follows as in lemma \ref{lem:vol} that $A_L$ is the approximate area of a cross section of $\Omega(\gamma)$. Indeed, for any $\gamma \in \G_L$,
 $vol(\Omega(\gamma)) = A_L length(\gamma) + O(Le^{-3L/2})$. So
\begin{eqnarray*}
\frac{\lambda(\Omega(\gamma))}{length(\gamma)} &=&
\frac{A_L}{vol(T_1(S))} + O(e^{-3L/2}).
\end{eqnarray*}
Thus
\begin{eqnarray*}
\frac{\lambda(\Omega_L)}{N} &=& \frac{\Sigma_{\gamma \in \G_L}
  \lambda(\Omega(\gamma))}{\Sigma_{\gamma \in \G_L} length(\gamma)}\\
&=&  \frac{\Sigma_{\gamma \in \G_L}
  [A_Llength(\gamma)/vol(T_1(S)) + O(Le^{-3L/2})]}{\Sigma_{\gamma \in \G_L} length(\gamma)}\\
 &=& A_L/vol(T_1(S)) + O(e^{-3L/2}).
\end{eqnarray*}
Let $\delta>0$. Choose $L$ larger if necessary so that $|f(x)-f(y)| < \delta$ if $distance(x,y) < 2\epsilon e^{-L/2}$. It follows that
\begin{eqnarray*}
| \omega_\gamma(f) - \frac{\lambda(\Omega(\gamma))}{length(\gamma)} \nu_\gamma(f)| \le \lambda(\Omega(\gamma)) \delta
\end{eqnarray*}
for any $\gamma \in \G_L$. But,

\begin{eqnarray*}
\Big|\frac{\lambda(\Omega_L)}{N} - \frac{\lambda(\Omega(\gamma))}{length(\gamma)}\Big|
 &\le&\Big|\frac{\lambda(\Omega_L)}{N} -  \frac{A_L}{vol(T_1(S))}\Big| + \Big| \frac{A_L}{vol(T_1(S))} -
 \frac{\lambda(\Omega(\gamma))}{length(\gamma)}\Big|\\
&=&  O(e^{-3L/2}).
\end{eqnarray*}
So
\begin{eqnarray*}
|\omega_\gamma(f) - \frac{\lambda(\Omega_L)}{N} \nu_\gamma(f)| &\le& \Big|\omega_\gamma(f) -\frac{\lambda(\Omega(\gamma))}{length(\gamma)} \nu_\gamma(f)\Big| + \Big| \frac{\lambda(\Omega(\gamma))}{length(\gamma)} \nu_\gamma(f) - \frac{\lambda(\Omega_L)}{N} \nu_\gamma(f)\Big|\\
&\le& \lambda(\Omega(\gamma)) \delta +  O(\nu_\gamma(f)e^{-3L/2}).
\end{eqnarray*}

Thus:

\begin{eqnarray*}
|\omega_L(f)-\nu_L(f)| &=& \big|\Sigma_{\gamma \in \G_L} \frac{1}{\lambda(\Omega_L)} \omega_\gamma(f)-\frac{1}{N}\nu_\gamma(f)\big|\\
                       &\le& \frac{1}{\lambda(\Omega_L)} \Sigma_{\gamma \in \G_L} |\omega_\gamma(f)-\frac{\lambda(\Omega_L)}{N}\nu_\gamma(f)|\\
                       &\le& \delta + O(\nu_L(\gamma)e^{-L/2}).
\end{eqnarray*}
Since $\delta>0$ is arbitrary, $|\omega_L(f)-\nu_L(f)| \to 0$ as $L$ tends to infinity as claimed. Since $f$ is arbitrary and $\omega_L$ converges to $\lambda$, $\nu_L$ also converges to $\lambda$.

\end{proof}

\section{Constructing Isometries in $\Gamma$}\label{section:isometry}

%the only thing left with this section is the signed distance statements of the corollary 2d.

Define
\begin{eqnarray*}
\sB_\epsilon=\{Mat[a,b,c,d]\in PSL_2(\R) | F(a,b,c) \in \{-\epsilon,\epsilon)\}
\end{eqnarray*}
where $F$ is as defined in corollary \ref{cor:perp}.

{\bf Motivating Problem}: given a vector $w \in T_1(\H^2)$, and $L,\epsilon>0$, find if possible an isometry $g \in \Gamma$ such that $tr.length(g)$ (the translation length of $g$) is in $(L-\epsilon,L+\epsilon)$ and the axis of $g$ is close to $axis(w)$, the geodesic in $\H^2$ tangent to $w$.

If there exists $g \in \Gamma$ such that $g wG_{-L/2} \in wG_{L/2} \sB_\epsilon$ then let $g_L(w)=g$. If $\epsilon$ is small enough (depending on the minimum injectivity radius of $S$) then $g_L(w)$ is unique if it exists. We always assume this is the case in what follows. Let $w_t = w R_{-\pi/2} G_t R_{\pi/2}$. $w_t$ is obtained by pushing $w$ to its right a distance $t$ along a geodesic. 

%$Comm(\Gamma)$, the commensurator of $\Gamma$ is the set of all $\gamma \in Isom(\H^2)$ such that $\gamma \Gamma \gamma^{-1} \cap \Gamma$ is a lattice.

\begin{cor}\label{cor:existence}
Assume $Comm(\Gamma)<Isom^+(\H^2)$. Then there exists a constant $L_0>0$ such that for any $L>L_0$ and any $w \in T_1(S)$ there exists a time $t$ with $0 \le t \le L_0e^{-L/2}$ such that $g_L(w_t)$ exists.
\end{cor}

\begin{proof}
Let $f: T_1(S) \times T_1(S) \to [0,1]$ be a continuous nonzero function with support contained in 
%\begin{eqnarray*}
%\big\{ (\Gamma g_1, \Gamma g_2) \in (\Gamma \backslash PSL_2(\R))\times  (\Gamma \backslash PSL_2(\R))| \Gamma g_2 \in \Gamma g_1 \sB_\epsilon\big\}.
%\end{eqnarray*}
\begin{eqnarray*}
\Big\{ (\Gamma g_1, \Gamma g_2) \, \big| \, \Gamma g_2 \in \Gamma g_1 \sB_\epsilon  \Big\}.
\end{eqnarray*}
By theorem \ref{thm:equidistribution} there exists a $L_0$ such that for all $T, L>L_0$ and $w \in T_1(\H^2)$ the support of $\pi_*(\mu_{w,T,L})$ intersects the support of $f$. This implies that for some $t$ with $0\le t <L_0e^{-L/2}$ (depending on $L>L_0$ and $w$), $w_{t}G_{-L/2} \in w_{t}G_{L/2} \sB_\epsilon$ and thus $g_L(w_t)$ exists. 
\end{proof}

%Note: if $g_L(w)$ exists then it is conjugate to $G_LMat[a,b,c,d]$ for some $Mat[a,b,c,d]\in \sB_\epsilon$.

The next result quantifies how close $axis(g_L(w))$ is to $axis(w)$ from the perspective of an observer on a geodesic $\tgamma$ parallel to $axis(w)$. It is proven in section \ref{section:hexagon H} as a corollary to theorem \ref{thm:main}.
 
\begin{figure}[htb] 
\begin{center}
 \ \psfig{file=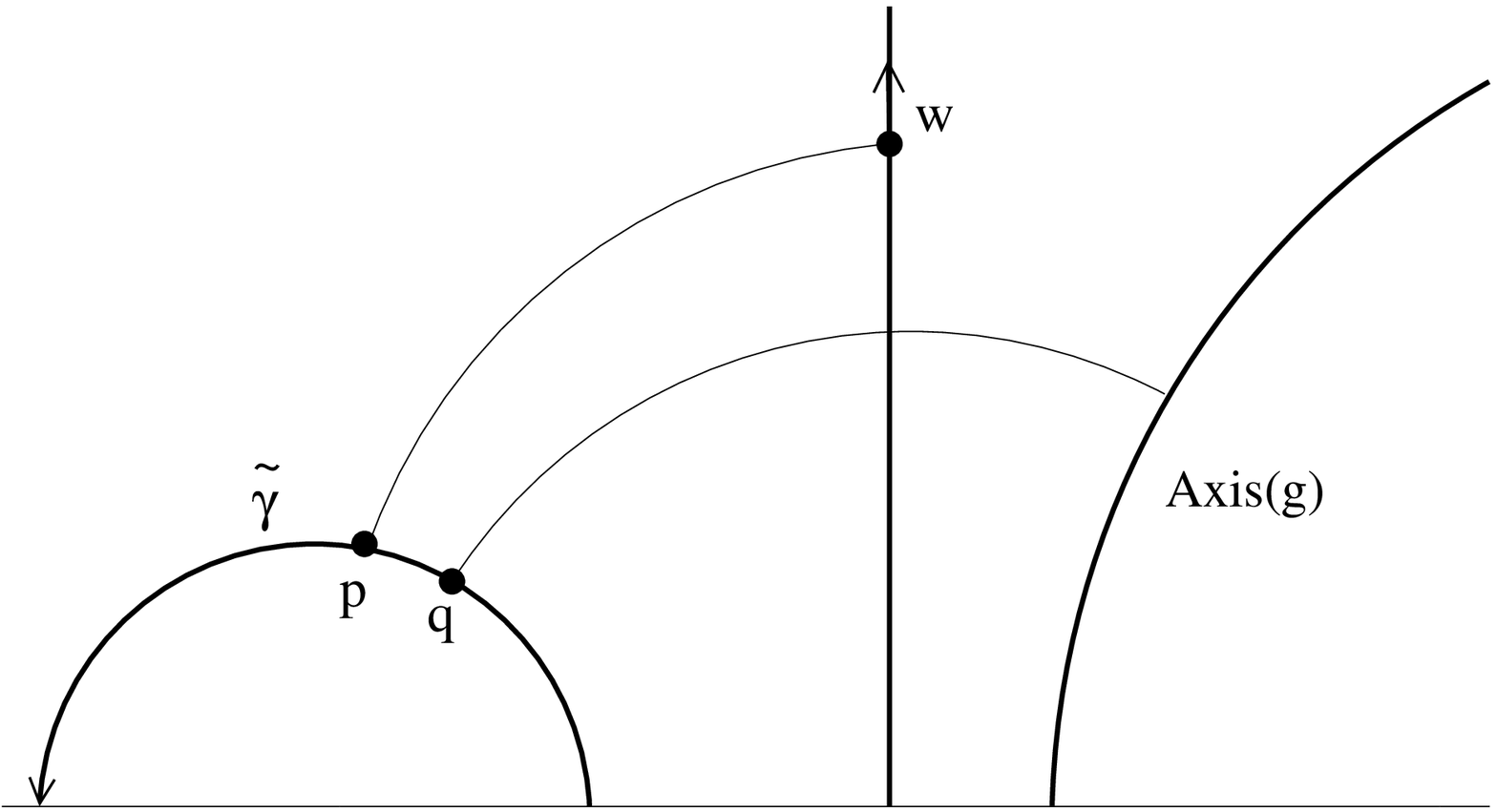,height=2in,width=3.7in}
 \caption{$\tgamma$, $w$ and $axis(g)$ in the upperhalf plane model.}
 \label{fig:hexagonH2}
 \end{center}
 \end{figure}

\begin{cor}\label{cor:2d}
Let $F=(F_1,F_2,F_3):(0,\infty) \times \R^2 \to \R^3$ be defined as in corollary \ref{cor:perp}. Suppose $w \in T_1(\H^2)$, $g:=g_L(w)$ exists and $gwG_{-L/2}=wG_{L/2}Mat[a,b,c,d]$ for some $a,b,c,d \in \R$ with $|ln(a)|,|b|,|c|,|d-1|<B$. Let $axis(w)$ be oriented consistently with $w$. Let $\tgamma$ be an oriented geodesic in $\H^2$ satisfying
\begin{itemize}
\item $\tgamma$ is on the left side of $axis(w)$ and
\item $axis(w)$ is on the right side of $\tgamma_1$.
\end{itemize}
See figure \ref{fig:hexagonH2}. Also assume that the shortest path from $\tgamma$ to $axis(w)$
\begin{itemize}
\item contains the basepoint of $w$,
\item intersects $\tgamma$ in a point $p$ and
\item has length $M=M(L)$ where either $M=m_1e^{-m_2L}$ for some $0\le m_2 <
 1/2$ or $M=m_3 + m_4L>0$ for some $m_4\ge 0$.
\end{itemize}
 Then the shortest path between $\tgamma$ and $axis(g)$ has length 
\begin{eqnarray*}
M+F_1(a,b,c)e^{-L/2} + O((B+1)^2\coth(M)e^{-L})
\end{eqnarray*}
and intersects $\tgamma$ in a
 point $q$. The signed distance from $p$ to $q$ equals
\begin{eqnarray*}
F_2(a,b,c)e^{-L/2}/\sinh(M) + O((B+1)^3 e^{-3L/2}/\sinh^3(M)).
\end{eqnarray*}
 The sign is positive depending on
 whether $q$ is before or after $p$ with respect to the orientation on
 $\tgamma$. The translation length of $g$ equals 
\begin{eqnarray*}
L + F_3(a,b,c) + O((B+1)e^{-L}).
\end{eqnarray*}
% All asymptotics here are taken as $L$ tends to infinity.
%If, instead $M= m_1L + m_2$ for some $m_1>0$ then the length of the shortest path between $\tgamma$ and $axis(g)$ has length asymptotic 
\end{cor}

\section{Constructing Pants}\label{section:pants}

%I'm pretty much happy with this whole section. pants16.tex.

Assume as given: positive numbers $r_1,r_2,r_3, L, \epsilon >0$, a closed oriented geodesic $\gamma_1$ in $S$ of length $l_1 \in (r_1 L -\epsilon, r_1L + \epsilon)$ and a point
$p$ on $\gamma_1$. 

{\bf Motivating Problem}: assuming $L$ is sufficiently large, construct an immersion of a pair of pants $P$ into $S$ such that
$\gamma_1$ is the image of one of the boundary components and the others $\gamma_2$, $\gamma_3$ have length in $r_2L + (-\epsilon,\epsilon)$, $r_3L + (-\epsilon,\epsilon)$ respectively. Also, the shortest path from $\gamma_1$ to $\gamma_2$ should have one endpoint close to $p$ and the orientation induced on $\gamma_1$ by $P$ should be its given orientation.

We will use the construction from the previous section to solve this problem under additional mild hypotheses. But first we need a lemma.

\begin{lem}\label{lem:pants}
Let $g_1,g_2 \in PSL_2(\R)$ be two hyperbolic isometries with disjoint
axes. For $i=1,2$ let  $l_i =
tr.length(g_i)$. Suppose there exists $l_3 \ge 0 $ satisfying
\begin{eqnarray*}
\cosh(l_3/2) = \sinh(l_1/2)\sinh(l_2/2)\cosh(M) - \cosh(l_1/2)\cosh(l_2/2)
\end{eqnarray*}
where $M$ is the length of the shortest path between $axis(g_1)$ and $axis(g_2)$. Let $C$ denote the convex hull of the limit set of $<g_1,g_2>$
(the group generated by $g_1,g_2$). Then the quotient $C/<g_1,g_2>$ is a hyperbolic
pair of pants $P$ and $l_1,l_2,l_3$ are the lengths of the boundary components.
\end{lem}
In the sequel, if $g_1,g_2 \in \Gamma$ satisfy the hypotheses above, we will say that $\{g_1,g_2\}$ determines $P$ and the immersion $j:P \to S$ induced by inclusion $<g_1,g_2> < \Gamma$. 

\begin{figure}[htb] 
\begin{center}
 \ \psfig{file=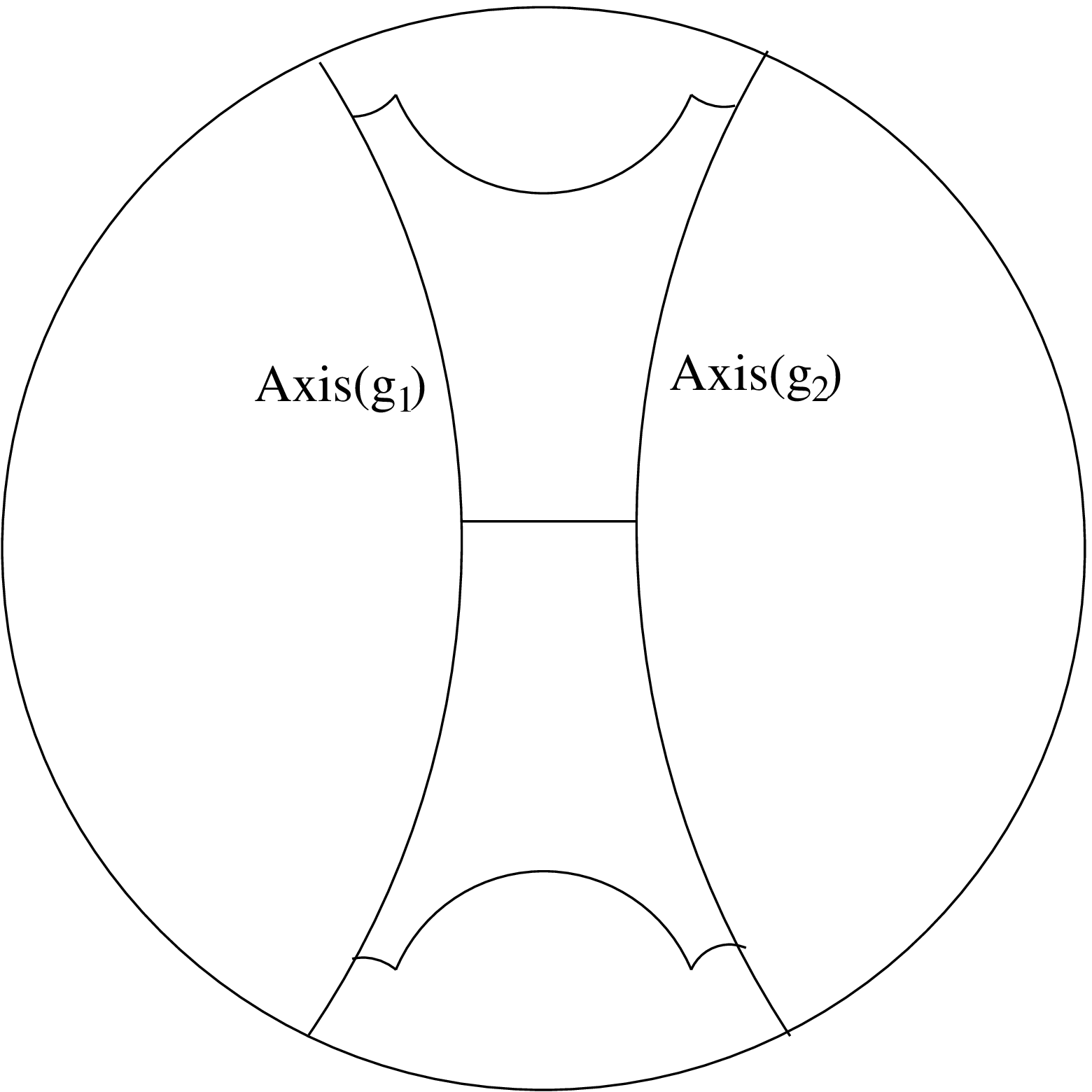,height=3in,width=3in}
 \caption{The hexagons $\sH_1$,$\sH_2$.}
 \label{fig:hexagons}
 \end{center}
 \end{figure}

\begin{proof}
There are two distinct right angled hexagons $\sH_1, \sH_2$ satisfying:
\begin{itemize}
\item $\sH_1, \sH_2$ are bounded by $axis(g_1)$, $axis(g_2)$ and the common perpendicular between $axis(g_1)$ and $axis(g_2)$
\item the sides of $\sH_1$ and $\sH_2$ contained in $axis(g_1)$ and $axis(g_2)$ each have length $tr.length(l_1)/2$ and $tr.length(l_2)/2$ respectively.
\end{itemize}
See figure \ref{fig:hexagons}. $\sF:=\sH_1 \cup \sH_2$ is a fundamental domain for the action of $<g_1,g_2>$ on the convex hull of its limit set. For $i=1,2$, the side of $\sF$ contained in $axis(g_i)$ has length $tr.length(g_i)$ and maps onto a boundary component of the quotient. Also the side of $\sH_i$ opposite the common perpendicular between $axis(g_1)$ and $axis(g_2)$ has length $l_3/2$. The above equation follows from the law of cosines (\cite{Rat}, \cite{Fen}).
\end{proof}

The next lemma is a short calculation. 
\begin{lem}\label{lem:pants estimates}
Assume the hypotheses of the previous lemma. In addition, suppose for $i=1,2$, 
\begin{eqnarray*}
l_i = r_iL + \rho_i
\end{eqnarray*}
where $r_1,r_2>0$ and $\rho_1,\rho_2 \in (-\epsilon,\epsilon)$. Let
 $r_3>0$ and $x\in \R$. Define $M=M(L)$ by:
\begin{displaymath}
M(L)=\left\{\begin{array}{ll}
    2e^{(-l_1 -r_2L+ r_3L+x)/4} & \textnormal{ if $r_1,r_2,r_3$ are the
    sidelengths of a Euclidean triangle}\\
    \arccosh(2e^{x/2}+1) & \textnormal{ if $r_1+r_2=r_3$}\\
    (-l_1-r_2L + r_3L + x)/2 + \ln(4) & \textnormal{ if $r_3 > r_1 + r_2$}
\end{array}\right.
\end{displaymath}
Then
\begin{displaymath}
l_3= r_3L + \rho_2 + x +  O(e^{-r_1L} + e^{-r_2L}).
\end{displaymath}

% Suppose $r_1,r_2,r_3$ are the sidelengths of a nondegenerate Euclidean triangle. If
%\begin{eqnarray*}
%M = 2e^{-l_1/4 -r_2L/4 + r_3L/4 + x/4}
%\end{eqnarray*}
%for some $x \in (-\epsilon,\epsilon)$ then
%\begin{eqnarray*}
%l_3 &=& r_3L + \rho_2 + x\\
%&& + O(e^{(-3r_1+r_2-r_3)L/2} + e^{(r_1-3r_2-r_3)L/2} + e^{(-r_1-r_2-r_3)L/2} + e^{-r_2L}+e^{-r_1L}).
%\end{eqnarray*}
%On the other hand, suppose $r_1 + r_2 < r_3$. If
%\begin{eqnarray*}
%M = (-l_1-r_2L + r_3L)/2 + \ln(4) + x/2
%\end{eqnarray*}
%for some $x \in (-\epsilon,\epsilon)$ then
%\begin{eqnarray*}
%l_3 &=& r_3L + \rho_2 + x + O(e^{-l_1} + e^{- l_2}).
%\end{eqnarray*}
\end{lem}

Back to the problem at hand: let $\tgamma_1, \tp$ denote lifts of $\gamma_1$ and
$p$ to the universal cover $\H^2$. Let $g_1 \in \Gamma$ be the element with axis $\tgamma_1$ and translation length equal to the length of $\gamma_1$. Let $M$ be defined as in the above lemma with $x=0$.

\begin{figure}[htb] 
\begin{center}
 \ \psfig{file=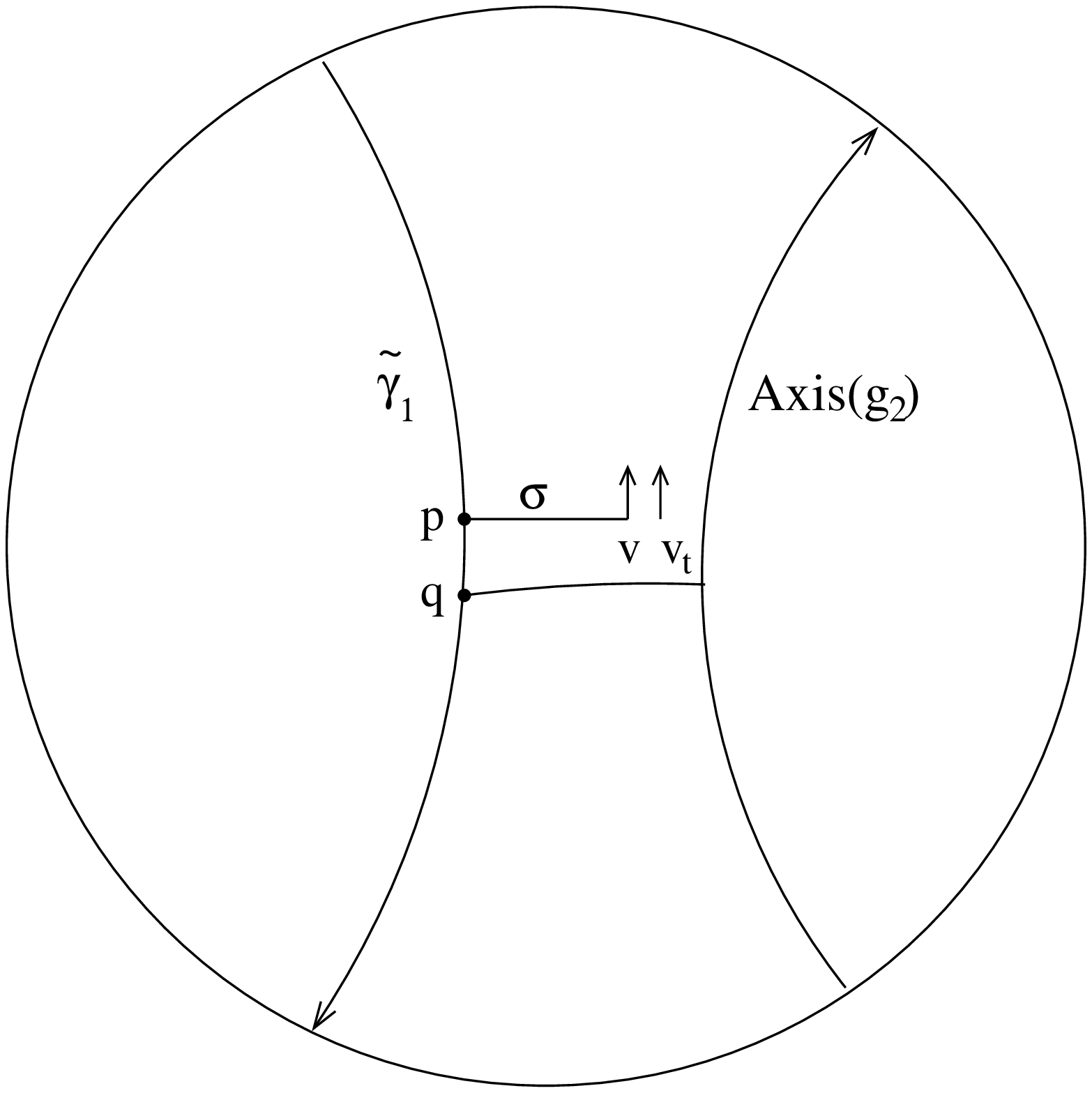,height=3in,width=3in}
\caption{}
 \label{fig:hexagons2}
 \end{center}
 \end{figure}

Let $\sigma$  be a segment of length $M$ orthogonal to $\tgamma_1$, with one endpoint at $\tp$ and on the left side of $\tgamma_1$. Let $v \in T_1(\H^2)$ be the unit vector based at the other endpoint of $\sigma$ that is orthogonal to $\sigma$ and such that $\sigma$ is on the left of $v$. Let $t$ be the smallest positive number such that $g_2:=g_{r_2L} (v_t)\in \Gamma$ exists where, as in the previous section, $v_t = vR_{-\pi/2} G_t R_{\pi/2}$. See figure \ref{fig:hexagons2}.

From now on, assume $Comm(\Gamma)<Isom^+(\H^2)$. By corollary \ref{cor:existence}, there exists an $L_0>0$ (depending only on $S$, $\epsilon$), so that if $L>L_0$ then $t \le L_0 e^{-r_2L/2}$. Let $\tq$ be the point of intersection between $\tgamma_1$ and the common perpendicular of $\tgamma_1$ and $\tgamma_2 = axis(g_2)$. Corollary \ref{cor:2d} implies:
\begin{cor}
If $L$ is sufficiently large then the distance between $axis(g_1)$ and $axis(g_2)$ equals $M+O(e^{-r_2L/2})$. The signed distance between $\tp$ and $\tq$ is $O(\epsilon e^{-r_2L/2}/\sinh(M))$. 
\end{cor}

By the previous lemma, $g_1$ and $g_2$ determine a pair of pants $P$ and an immersion $j: P \to S$ induced by inclusion $<g_1,g_2> < \Gamma$. Lemma \ref{lem:pants estimates} now implies

\begin{cor}
If $Comm(\Gamma)<Isom^+(\H^2)$ then
\begin{itemize}
 \item for all sufficiently large $L$,
\item  for any closed oriented geodesic $\gamma_1$ in $S$ with length in $r_1L + (-\epsilon,+\epsilon)$,
 \item and for any point $p$ in $\gamma_1$
\end{itemize}
there exists an immersion $j:P\to S$ of a pair of pants such that if $c_1,c_2,c_3$ are the boundary components of $P$ then
\begin{itemize}
\item $length(c_i) \in r_iL + (-\epsilon,\epsilon)$ (for $i=1,2,3$),
\item $c_1$ maps onto $\gamma_1$ in an orientation-preserving way,
\item if $q$ is the intersection point between $c_1$ and the shortest arc in $P$ between $c_1$ and $c_2$ then the distance between $j(q)$ and $p$ is $O(\epsilon e^{-r_2L}/\sinh(M))$.
\end{itemize}
\end{cor}

%\begin{remark}
%It may be possible to remove some of the hypotheses on $r_1,r_2,r_3$ above although we have not tried.
%\end{remark}

The next corollary yields a weak form of the Ehrenpreis conjecture. For this paper, a pants decomposition of a hyperbolic surface $\Sigma$ is a collection $C$ of disjoint geodesics such that every component of the complement $\Sigma - \cup_{\gamma \in C} \gamma$ is a pair of pants with finite area. For every geodesic $\gamma \in C$, the twist parameter at $\gamma$ is defined as follows. Let $P_1, P_2$ be two different components of the complement $\Sigma - \cup_{\gamma \in C} \gamma$ whose closures contain $\gamma$. For $i=1,2$ let $p_i \in \gamma$ be the endpoint of the shortest arc between $\gamma$ and some other boundary component of $P_i$. Then $twist(\gamma)$ equals the distance between $p_1$ and $p_2$ (mod $length(\gamma)/2$).

\begin{cor}\label{cor:ehrenpreis}
For any closed surface $S=\H^2/\Gamma$ with $\Gamma, Comm(\Gamma) < Isom^+(\H^2)$, any $\epsilon>0$ and for all sufficiently large $L$, there exists a locally-isometric covering map $j: {\tilde S} \to S$ where ${\tilde S}$ has a pants decomposition in which every geodesic of the decomposition has length in $(L-\epsilon,L+\epsilon)$ and every twist parameter is bounded by $\epsilon e^{-L/4}$.
\end{cor}
This corollary follows from the previous corollary by setting $r_1=r_2=r_3=1$ and successively gluing together immersions of pants. We cannot guarantee that ${\tilde S}$ can be chosen to be closed; if we could, the Ehrenpreis conjecture would follow shortly. Using standard calculations, it can be shown that there exists a $(1+10\epsilon)$-quasiconformal homeomorphism $h: {\tilde S} \to {\hat S}$ where ${\hat S}$ is a surface with a pants decomposition in which every geodesic has length exactly $L$ and every twist parameter is exactly zero.

%\begin{remark}
%More generally, it is possible to require ${\tilde S}$ to have any pre-chosen twist parameters whatsoever, allowing for an error on the order of $\epsilon e^{-L/4}$. As in the previous remark, there exists a $(1+10\epsilon)$-quasiconformal homeomorphism such a ${\tilde S}$ to a surface with a pants decomposition in which every geodesic has length exactly $L$ and every twist parameter is exactly as as pre-specified.%
%\end{remark}

\section{Counting}\label{section:counting}

The goal of this section is to prove theorem \ref{thm:counting}, corollary \ref{cor:all} and theorem \ref{thm:twist}.

%----------------------------
\begin{proof}(of theorem \ref{thm:counting})
 Let $L>0$. Let $\gamma \in \G_{r_1L}$. Let $\tgamma$ be a lift of
 $\gamma$ to $\H^2$. Let $g \in \Gamma$ be the element with axis
 $\tgamma$ and translation length equal to $l_1:=length(\gamma)$. Let $\P_L(\gamma)=\P_L(r_1,r_2,r_3;\gamma)$.

Suppose $(j:P \to S) \in \P_L(\gamma)$. Then there exists
an element $h \in \Gamma$ such that $\{g,h\}$ determines $P$ in the
sense of lemma \ref{lem:pants} and $tr.length(h) \in r_2L + (-\epsilon,\epsilon)$. Since the third boundary component has length $l_3 \in r_3L + (-\epsilon, \epsilon)$ lemma \ref{lem:pants} implies that the shortest distance $M$ between $axis(h)$ and $axis(g)$ is in the interval $(M_-,M_+)$ where 
\begin{eqnarray*}
\cosh(M_{\pm}) &=& \frac{\cosh((r_3L\pm \epsilon)/2) + \cosh(l_1/2)\cosh((r_2L \mp \epsilon)/2)}{\sinh(l_1/2)\sinh((r_2L\mp \epsilon)/2)}.
\end{eqnarray*}
To interpolate between the two values, it is convenient to define $M'_x$ (for $x \in [-\epsilon,\epsilon]$) by 
\begin{eqnarray*}
\cosh(M'_x) &=& \frac{\cosh((r_3L+x)/2) + \cosh(l_1/2)\cosh((r_2L-x)/2)}{\sinh(l_1/2)\sinh((r_2L-x)/2)}\\
           &=& 2e^{(-l_1 -r_2L + r_3L)/2}e^x +1 + O(E_2).
\end{eqnarray*}
where the error term, $E_2$ is on the order of
\begin{eqnarray*}
e^{(-r_1-r_2-r_3)L/2} + e^{-r_1L} + e^{-r_2L} + e^{(-3r_1-r_2+r_3)L/2}+ e^{(-r_1-3r_2+r_3)L/2}  .
\end{eqnarray*}
Since the error term above is asymptotically negligible compared with the main term, we will ignore it in the calculations below. Define
 \begin{eqnarray*}
\cosh(M_x) &=& 2e^{(-l_1 -r_2L + r_3L)/2}e^x +1.
\end{eqnarray*}
So $M_-\sim M_{-\epsilon}$ and $M_+ \sim M_{\epsilon}$. Now suppose that the shortest distance between $axis(h)$ and $axis(g)$ equals $M_x$ for some $x \in [-\epsilon,\epsilon]$. Let $l_2 = tr.length(h)=L+\rho_2$ for some fixed $\rho_2 \in [-\epsilon,\epsilon]$. Then
 \begin{eqnarray*}
2e^{(-l_1 -r_2L + r_3L)/2}e^x + 1 &\sim & 2e^{(-l_1 -l_2 + l_3)/2} + 1.
\end{eqnarray*}
So $-r_2L + r_3L + 2x \sim -l_2 + l_3$. Since $l_3 \in r_3L + (-\epsilon,\epsilon)$, we must have
\begin{eqnarray*}
l_2 - r_2L \in -2x+ (-\epsilon,\epsilon).
\end{eqnarray*}
Since $\{g,h\}$ determines a pants immersion in $\P(\gamma)$, we must also have
\begin{eqnarray*}
l_2 -r_2L &\in& (-2x-\epsilon,-2x + \epsilon) \cap (-\epsilon,\epsilon)\\
          &=&  \big(\max(-\epsilon-2x,-\epsilon),\min(\epsilon-2x,\epsilon)\big).
\end{eqnarray*}
So define:
\begin{eqnarray*}
X(x)=  (-1,1) \times (-1,1) \times \big(\max(-\epsilon-2x,-\epsilon),\min(\epsilon-2x,\epsilon)\big)
\end{eqnarray*}
and
\begin{eqnarray*}
\sB(x) &=& \{Mat[a,b,c,d] \in PSL_2(\R): F(a,b,c) \in X(x)\}
\end{eqnarray*}
where $F$ is as defined in corollary \ref{cor:perp}. For $w \in T_1(\H^2)$ let $z(w)$ denote the signed distance from $basepoint(w)$ to $axis(g)$ where the sign is positive iff $w$ is on the left side of $axis(g)$. Let $x(w)$ be defined by $z(w)=M_{x(w)}$ when this is defined. For $w \in T_1(\H^2)$, if there exists $h \in \Gamma$ such that
\begin{eqnarray*}
 h wG_{-r_2L/2} \in wG_{r_2L/2} \sB(x(w))
\end{eqnarray*}
then define $h_w = h$. Let $V \subset T_1(\H^2)$ be the set of all vectors $w \in T_1(\H^2)$ satisfying:
\begin{itemize}
\item $x(w) \in [-\epsilon,\epsilon]$.
\item the shortest path from the basepoint of $w$ to $axis(g)$ is perpendicular to $w$ and is on the left side of $\tgamma$.
\end{itemize}
Let $W$ be the set of all vectors $w \in V$ such that $h_w$ exists. Ignoring the asymptotically neglible error terms above, we have shown that $\{g,h\}$ determines a pants immersion in $\P(\gamma)$ iff $h=h_w$ for some $w$ in $W$. For $h \in \Gamma$, define
\begin{eqnarray*}
W(h) = \{w \in W: h_w = h\}.
\end{eqnarray*}
Let $\pi: T_1(\H_2) \to T_1(S)$ be the quotient map. Let $w$ be an element of $\pi(W)$ chosen uniformly at random. Then
\begin{eqnarray*}
n|\P(\gamma)| &=& area(\pi(W)) \mE\Big(\frac{1}{area(W(h_w))}\Big).
\end{eqnarray*}
Here $n$ equals $1$ if $r_2 \ne r_3$ and $2$ otherwise. $\mE$ is expectation with respect to the law of $w$. We have abused notation by writing $area(W(h_w))$ for $area(W(h_{\tilde w}))$ where ${\tilde w} \in W$ is a lift of $w$. We will use theorem \ref{thm:equidistribution} to obtain asymptotic estimates for the right hand side of the above equation. 

Let $h=h_w$ for some $w \in W$. By definition of $\sB(x)$ and corollary \ref{cor:2d}, $W(h)$ is approximately rectangular. Its projection to $\H^2$ is, up to first order, a union of hypercycle segments $s_z$ for $z \in z_h + (-e^{-r_2L/2},e^{-r_2L/2})$ where 
\begin{itemize}
\item $z_h$ is the shortest distance from $axis(h)$ to $axis(g)$,
\item the distance from any point in $s_z$ to $axis(g)$ equals $z$,
\item the projection of $s_z$ to $axis(g)$ has length $\frac{2e^{-r_2L/2}}{\sinh(z)}$.
\end{itemize}
See figure \ref{fig:basin}.

\begin{figure}[htb] 
\begin{center}
 \ \psfig{file=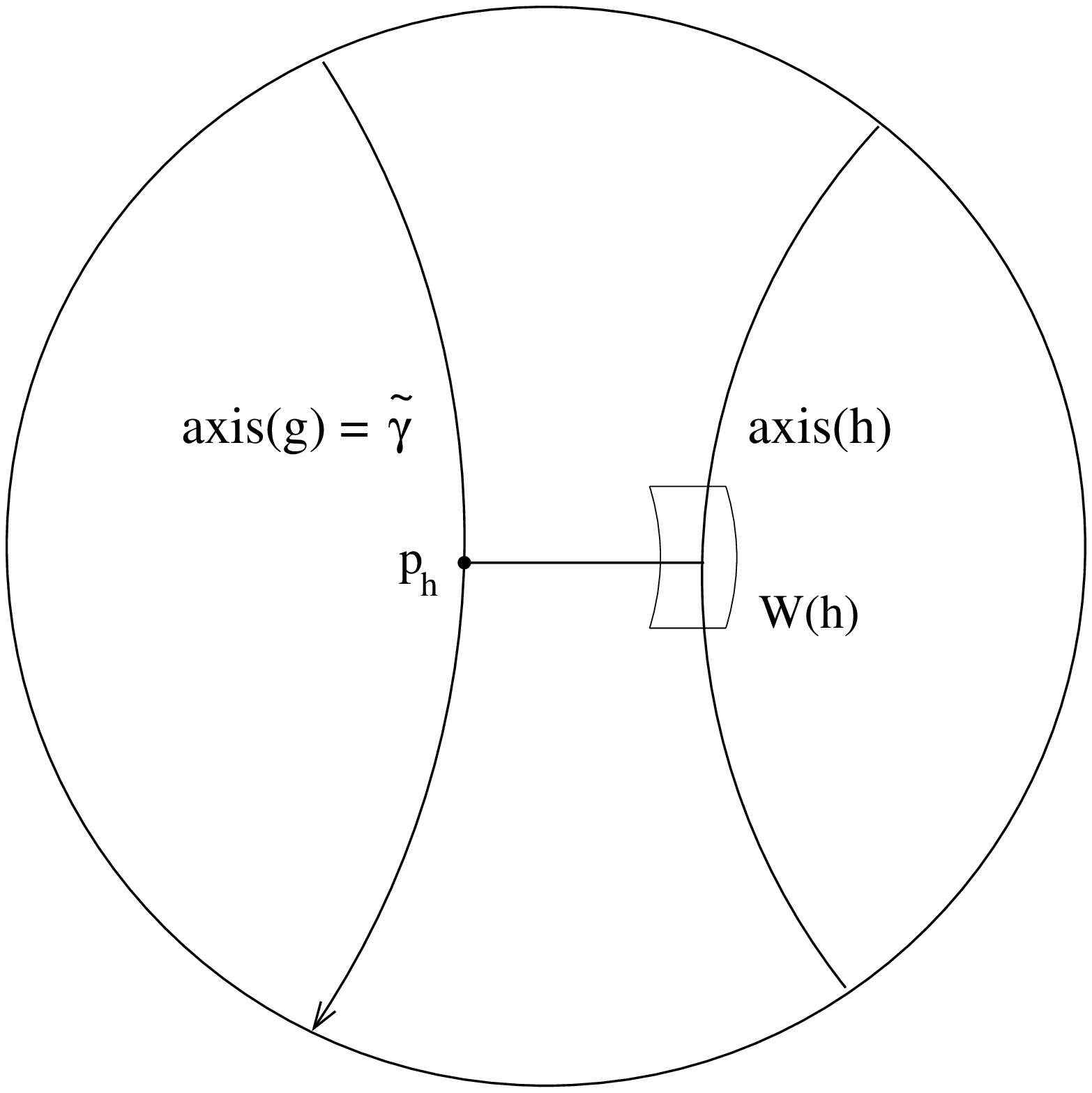,height=3in,width=3in}
 \caption{}
 \label{fig:basin}
 \end{center}
 \end{figure}

 The length of $s_z$ equals $2e^{-r_2L/2}\coth(z)$. So 
\begin{eqnarray*}
area(W(h)) &\sim& \int_{z_h - e^{-r_2L/2}}^{z_h + e^{-r_2L/2}} 2e^{-r_2L/2}\coth(z) dz\\
           &=& 2e^{-r_2L/2} \ln\Big(\frac{\sinh(z_h + e^{-r_2L/2})}{\sinh(z_h - e^{-r_2L/2})}\Big).
\end{eqnarray*}
For $z \in [M_-,M_+]$ define
\begin{eqnarray*}
A(z) =  2e^{-r_2L/2} \ln\Big(\frac{\sinh(z + e^{-r_2L/2})}{\sinh(z - e^{-r_2L/2})}\Big).
\end{eqnarray*}

Now if $v \in \pi(V)$ is chosen uniformly at random and $z(v):=z(\tv)$ where $\tv \in V$ is any lift of $v$ then the distribution of $z(v)$ is supported on $[M_-,M_+]$ with density equal to
\begin{eqnarray*}
\frac{length(\bgamma) \cosh(z)}{area(\pi(V))}.
\end{eqnarray*}

Similarly, if $w \in \pi(W)$ is chosen uniformly at random, then the distribution of $z(w)$ is supported in $[M_-,M_+]$ with density asymptotic to
\begin{eqnarray*}
\frac{length(\bgamma)}{area(\pi(W))} \cosh(z) \frac{vol(\sB(x))}{vol(T_1(S))}
\end{eqnarray*}
where $M_x = z$. This uses theorem \ref{thm:equidistribution}. So
\begin{eqnarray*}
 n|\P(\gamma)| &=& area(\pi(W)) \mE\Big(\frac{1}{area(W(h_w))}\Big)\\
               &\sim&  \int_{M_-}^{M_+} \, \frac{length(\bgamma)}{A(z)} \cosh(z) \frac{vol(\sB(x))}{vol(T_1(S))} dz.
\end{eqnarray*}
Since
\begin{eqnarray*}
\cosh(z) &=& 2e^{(-l_1 -r_2L + r_3L)/2}e^x + 1,
\end{eqnarray*}
\begin{eqnarray*}
\sinh(z) dz &=& 2e^{(-l_1 -r_2L + r_3L)/2}e^x dx.
\end{eqnarray*}
So,
\begin{eqnarray*}
 n|\P(\gamma)| &\sim& \frac{2e^{(-l_1 -r_2L + r_3L)/2}length(\bgamma)}{vol(T_1(S))} \int_{-\epsilon}^{\epsilon} \, \frac{1}{A(M_x)} \coth(M_x) e^x vol(\sB(x)) dx.
\end{eqnarray*}
A short calculation shows:
\begin{eqnarray*}
A(M_x) \sim 4e^{-r_2L} \coth(M_x).
\end{eqnarray*}
Hence
\begin{eqnarray*}
 n|\P(\gamma)| &\sim& \frac{e^{(-l_1 +r_2L + r_3L)/2}length(\bgamma)}{2vol(T_1(S))} \int_{-\epsilon}^{\epsilon} \, e^x vol(\sB(x)) dx.
\end{eqnarray*}

By lemma \ref{lem:vol},
\begin{displaymath}
vol(\sB(x))=\left\{\begin{array}{ll}
    4(e^{\epsilon-2x}-e^{-\epsilon}) & \textnormal{ if $x>0$}\\
    4(e^{\epsilon}-e^{-\epsilon-2x})  & \textnormal{ if $x<0$}.
\end{array}\right.
\end{displaymath}
So $vol(\sB_x) = 8e^{-x} \sinh(\epsilon - |x|)$. Hence
\begin{eqnarray*}
\int_{-\epsilon}^{\epsilon} \, e^x vol(\sB(x)) dx &=& \int_{-\epsilon}^{\epsilon} \, 8\sinh(\epsilon - |x|) dx\\
%&=& 8\big(-\cosh(0)+\cosh(\epsilon) + \cosh(\epsilon) -\cosh(0)\big)\\
%&=& 16( \cosh(\epsilon)-1)\\
&=& 8(e^{\epsilon/2}-e^{-\epsilon/2})^2.
\end{eqnarray*}
Hence
\begin{eqnarray*}
 n|\P(\gamma)| &\sim& \frac{4(e^{\epsilon/2} -e^{-\epsilon/2})^2 e^{(-l_1 +r_2L + r_3L)/2}length(\bgamma)}{vol(T_1(S))}.
\end{eqnarray*}

\end{proof}

%The function $f$ from the statement of the above theorem, is:
%\begin{eqnarray*}
%%\frac{n}{8\pi \epsilon^2} \int_{-\epsilon}^{\epsilon} e^{x/4} vol(\sB'_\epsilon(x)) dx
%\end{eqnarray*}
%in the case $r_1 + r_2 > r_3$.

\begin{proof}(of corollary \ref{cor:all})
Assume, for simplicity, that $r_1 < r_2 < r_3$. Then $\P(r_1,r_2,r_3)$ is the disjoint union of $\P(r_1,r_2,r_3;\gamma)$ over all $\gamma \in \G_{r_1L}$. Let $\G'_{r_1L} \subset \G_{r_1L}$ be the subset of geodesics that trivially cover their images. If $\gamma$ is a uniformly random element of $\G'_{r_1L}$, then the distribution of $length(\gamma)-r_1L$ has support in $[-\epsilon,\epsilon]$ and density asymptotic to
\begin{eqnarray*}
\frac{1}{e^\epsilon-e^{-\epsilon}} e^z
\end{eqnarray*}
Hence
%For $x<y$, let $\G_{r_1L}(x,y)$ denote the set of closed geodesics in $S$ which have length in $r_1L + (x,y)$. Let $\G'_{r_1L}(x,y) \subset \G_{r_1L}(x,y)$ be the subset of geodesics that trivially cover their images. It is well-known that
%\begin{eqnarray*}
%|\G'_{r_1L}(x,y)| \sim |\G_{r_1L}(x,y)| \sim  (e^\epsilon-e^{-\epsilon})e^{r_1L}/(r_1L).
%\end{eqnarray*}
%Hence
\begin{eqnarray*}
|\P(r_1,r_2,r_3)| &=& \Sigma_{\gamma \in \G_{r_1L}} |\P(r_1,r_2,r_3;\gamma)|\\
&\sim & \Sigma_{\gamma \in \G'_{r_1L}} |\P(r_1,r_2,r_3;\gamma)|\\
&\sim& |\G'_{r_1L}| \mE(\P(r_1,r_2,r_3;\gamma))
\end{eqnarray*}
where $\mE$ is expectation with respect to the uniform measure on $\gamma \in \G'_{r_1L}$. By the above,
\begin{eqnarray*}
\mE(\P(r_1,r_2,r_3;\gamma)) &\sim& \frac{4(e^{\epsilon/2}-e^{-\epsilon/2})^2 }{ vol(T_1(S)) }\frac{1}{e^\epsilon-e^{-\epsilon}} \int_{-\epsilon}^{\epsilon} (r_1L +z)e^{(-(r_1L+z) +r_2L + r_3L)/2}e^z \, dz\\
&=& \frac{4(e^{\epsilon/2}-e^{-\epsilon/2})^2 e^{(-r_1+r_2+r_3)L/2} }{(e^\epsilon-e^{-\epsilon}) vol(T_1(S)) } \int_{-\epsilon}^{\epsilon} (r_1L +z)e^{z/2} \, dz.
\end{eqnarray*}
But
\begin{eqnarray*}
\int_{-\epsilon}^{\epsilon} (r_1L +z)e^{z/2} \, dz &=& 2(r_1L-2+\epsilon)e^{\epsilon/2} - 2(r_1L-2-\epsilon)e^{-\epsilon/2}.
\end{eqnarray*}
Combined with $|\G'_{r_1L}| \sim (e^\epsilon-e^{-\epsilon}) e^{r_1L}/(r_1L)$ the last two equations yield the corollary.

\end{proof}

\begin{proof}(of theorem \ref{thm:twist})
We will only prove the theorem in the case that $r_1 + r_2 > r_3$. The other cases are handled similarly. Let $\gamma, \tgamma, g, M$ be as defined in the previous proof.

Let $H \subset \Gamma$ be the set of all elements $h$ such that $\{h,g\}$ determines a pants immersion $(j:P \to S) \in \P_L(r_1,r_2,r_3;\gamma)$. For each $h \in H$, let $p_h \in \tgamma$ be the point closest to $axis(h)$. Let $\mu_L$ be the counting measure supported on the points $\{p_h: h \in H\}$. We will show that appropriately normalized, $\mu_L$ converges to Lebesgue measure.

Let $\sO_L$ be any point in $\tgamma_L$. This choice and the orientation of $\tgamma_L$ determines an isometry between $\tgamma_L$ and the real line sending $\sO_L$ to the origin. Identify $\tgamma_L$ with the real line through this isometry so that $\mu_L$ can now be thought of as a measure on the real line. 

Truncate $\mu_L$ by setting $\mu'_L := \chi_{[-r,r]} \mu_L$ where $2r=length(\bgamma)$ and $\chi_{[-r,r]}$ is the characteristic function of $[-r,r]$. Next dilate to the scale of interest by setting $\mu''_L(E) := \mu'_L(s E)$ where $E$ is any Borel set in $\R$ and $s = \epsilon e^{-r_2L/2}/\sinh(M)$. Now $\mu''_L$ is supported on $[-r/s, r/s]$. Normalize by setting $\mu'''_L = \frac{1}{K} \mu''_L$ where $K$ is chosen so that $\mu'''(L)$ has total mass $2r/s$.

To prove the theorem, it is necessary and sufficient to show that $\mu'''_L$ converges to Lebesgue measure in the weak* topology regardless of the choice of origin $\sO_L$. But this is equivalent to stating that the measure $\nu_L$ defined by $\nu_L =\chi_{[-1,1]} \mu'''_L$ converges to Lebesgue measure on the interval $[-1,1]$. 

%**define $\omega_L$

Recall the definition of $W(h)$ from the proof of the previous theorem. Let $\sigma$ be the geodesic orthogonal to $\tgamma$ that contains $\sO_L$. Let $\omega_L$ be the discrete measure on $[-1,1]$ defined setting $\omega_L(\{p\})$ equal to the length of 
\begin{eqnarray*}
\bigcup [W(h) \cap \sigma]
\end{eqnarray*}
where the union is over all $h \in H$ such that the shortest path from $axis(h)$ to $\tgamma$ contains $sp$. We have abused notation here by identifying $W(h)\subset T_1(\H^2)$ with its projection to the plane $\H^2$. 
%**relationship between $\omega_L$ and $\mu'''_L$

Recall that, up to a negligible error, $W(h)$ is a ``rectangle'' bounded by two geodesics separated a distance $2s$ apart and two hypercycles separated a distance $2\epsilon e^{-r_2L/2}$ apart. So for $h\in H$ either 
\begin{itemize}
\item the length of $W(h)\cap \sigma$ is asymptotically negligible (for example if $W(h)\cap \sigma = \emptyset$),
\item the length of $W(h)\cap \sigma$ equals $2\epsilon e^{-r_2L/2} + O(\coth(M)e^{-r_2L})$ (by corollary \ref{cor:2d} and the definition of $W(h)$)
\item or the distance from $axis(h)$ to $\tgamma$, is within $2\epsilon e^{-r_2L/2}$ of $Me^{-2\epsilon/4}$ or $Me^{2\epsilon/4}$.
\end{itemize}
The last case contributes an asymptotically negligible amount to the calculations. The support of $\omega_L$ is essentially the same as the support of $\nu_L$ because $\nu_L(\{p\}) \ne 0$ iff there exists an element $h\in H$ such that the shortest path from $axis(h)$ to $axis(g)$ contains $sp$. But this occurs iff $W(h) \cap \sigma \ne \emptyset$ since $W(h)$ is, up to a negligible error, bounded by two geodesics orthogonal to $\tgamma$ separated a distance $2s$ apart with common perpendicular midpoint at $sp$. So if $\omega_L$ is normalized to have total mass $2$ then any weak* limit point of $\omega_L$ is a weak* limit point of $\nu_L$ and vice-versa. So it suffices to show that $\omega_L$ limits on Lebesgue measure.

%**$\omega_L$, $y_v$ and equidistribution.

Let $v \in W \cap \sigma$ be chosen uniformly at random. Let $y_v \in \R$ be such that the distance from $\sO$ to $p$ equals $y_v s$ where $p$ is the closest point on $\tgamma$ to $axis(g_v)$. So the distribution of $y_v$ is the same as the measure $\omega_L$ (after normalizing to have total mass $1$). By corollary \ref{cor:2d}, $y_v = F_2(a,b,c)\sinh(M)/\sinh(z_v)$ where $z_v$ is the distance from the basepoint of $v$ to $\tgamma$ and $(a,b,c) \in \R^3$ is such that
\begin{eqnarray*}
g_v v G_{-L/2} = v G_{L/2} Mat[a,b,c,d].
\end{eqnarray*}
If $z_v=Me^{x_v/4}$ then $|x_v|<\epsilon$ and $y_v$ converges to $F_2(a,b,c) e^{-x_v/4}$. Asymptotically then, $y_v$ only depends on $a,b,c$ and $x_v$. So it follows from the equidistribution theorem that the distribution of $y_v$, and therefore $\omega_L$, converges to a limiting measure that is independent of origin $\sO$. Independence implies translation invariance. Hence this measure must be Lebesgue measure.

\end{proof}

%------------------

\section{Clotheslines}\label{section:clotheslines}

%we should go over this section again pants16.

The goal of this section is to prove theorem \ref{thm:clotheslines}. We will need the following estimate.

\begin{lem}\label{lem:x}
Suppose $w \in T_1(\H^2)$, $g \in \Gamma$ exists such that $gwG_{-L/2}=wG_{L/2}Mat[a,b,c,d]$ for some $a,b,c,d \in \R$. Let $axis(w) \subset \H^2$ be the geodesic tangent to $w$. Suppose $\gamma \subset \H^2$ is a geodesic such that the shortest path $\beta$ from $\gamma$ to $axis(w)$ satisfies
\begin{itemize}
\item $\beta$ has length $M=m_1e^{-m_2L}$ (for some $m_1,m_2\ge 0$, $m_2 < 1/2$),
\item the distance from $\beta \cap axis(w)$ to $basepoint(w)$ equals $K=k_1+ k_2 L$ for some $k_1 \in \R$, $k_2 \ge 0$. 
\end{itemize}
Then the shortest path from $\gamma$ to $axis(g)$ has length $M + O(\frac{|c|+|b|+1}{a}e^{k_1+(k_2-1/2)L}) $.
\end{lem}

\begin{proof}
After conjugating, we may assume that $w=(e^{L/2}i,e^{L/2}i)$ in the upperhalf plane model. Let $J=e^{\pm K+L/2 }$. The lemma now follows from theorem \ref{thm:main} below.
\end{proof}

\begin{defn}
Given $\delta >0$, a set $X \subset T_1(S)$ is called
$\delta$-dense if for every $v \in T_1(S)$ there exists a geodesic segment $\beta$ such that
\begin{itemize}
\item $length(\beta) < \delta$,
\item one endpoint of $\beta$ is at $basepoint(v)$ the other is in $basepoint(x)$ for some $x \in X$,
\item $\beta$ is perpendicular to $v$ and $v$ parallel transported across $\beta$ makes an angle $\le \delta$ with $x$.
\end{itemize}
% at most $\delta$ such that exists an $x \in
%X$ such that $d(v,x) < \delta$. Here $d(\cdot,\cdot)$ denotes the distance function in $T_1(S)$.
\end{defn}
For every  oriented closed geodesic $\gamma$ in either $S$ or $\H^2$, identify $\gamma$
with the set of unit vectors $v$ tangent to $\gamma$ and
oriented consistently with $\gamma$.

% Recall that $\G_L$ is the set of oriented closed geodesics in $S$ of length in $(L-\epsilon,L+\epsilon)$ and $\P_L(1,1,1)$ is the set of immersion of pants $j:P\to S$ such that all boundary components of $P$ map to geodesics in $\G_L$.

\begin{lem}
There exists a $\delta > 0$ (depending only on $S$) such that if $L$
is sufficiently large, $\gamma_1 \in \G_L$ is $\delta$-dense and
$\gamma_2 \in \G_L$ is arbitrary then there
exists an $L$-clothesline of length 2 from $\gamma_1$ to $\gamma_2$.  
\end{lem}

\begin{proof}
Let $\gamma_1,\gamma_2 \in \G_L$ be given. Suppose $\gamma_1$ is
$\delta$-dense for some $\pi/2>\delta>0$ with $\cosh(2\delta)\cos(\delta) >
1$. Assume $L$ is large enough so that if
$C=4e^{-L/4+\epsilon/4}$ then $C < \delta$ and $\cosh(C) < \cosh(2\delta)\cos(\delta)$.

Claim: there exists a geodesic segment $\beta$ perpendicular to $\gamma_1$ and $\gamma_2$ at
its endpoints so that 
\begin{itemize}
\item $4\delta\ge length(\beta)> C$ and
\item $\beta$ is to the left side of both $\gamma_1$ and $\gamma_2$.
\end{itemize}

Proof: Let $\beta_1$ be any geodesic segment of length $3\delta$ with one endpoint making a right angle with $\gamma_2$ such that $\beta_1$ is on the left side of $\gamma_2$. See figure \ref{fig:dense}. Let $u$ be the unit vector at that endpoint that points in the direction of $\gamma_2$. Let $v$ equal $u$ parallel transported across $\beta_1$. Since $\gamma_1$ is $\delta$-dense, there exists a geodesic segment $\beta_2$ of length at most $\delta$ with one endpoint making a right angle with $v$ and the other in $\gamma_1$ so that $-v$ parallel transported across $\beta_2$ makes an angle at most $\delta$ with a vector in $\gamma_1$. The concatenation of $\beta_1$ with $\beta_2$ gives a geodesic segment $\beta_3$ of length between $2\delta$ and $4\delta$ such that $\beta_3$ makes a right angle with $\gamma_2$ and an angle $\pi/2-\alpha$ with $\gamma_1$ where $|\alpha|<\delta$.

\begin{figure}[htb] 
\begin{center}
 \ \psfig{file=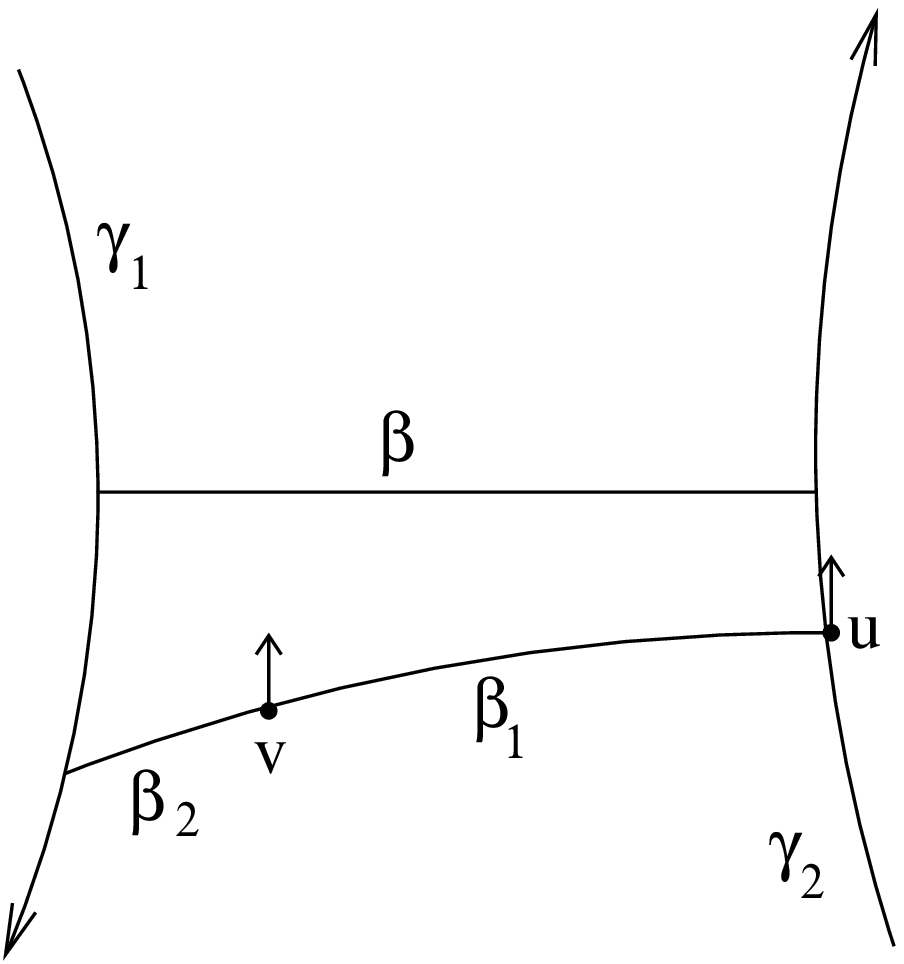,height=2.2in,width=2.2in}
 \caption{$\gamma_1,\gamma_2, \beta$.}
 \label{fig:dense}
 \end{center}
 \end{figure}

Let $\beta$ be the perpendicular segment from $\gamma_1$ to $\gamma_2$
that is homotopic to $\beta_3$ through homotopies that keep the
endpoints of $\beta_3$ in $\gamma_1$ and $\gamma_2$. Then $\beta$ and
$\beta_3$ are opposite sides of a 4-gon with three right angles and
one angle equal to $\pi-|\alpha|$. By the
trirectangle identities \cite{Bus}, 
\begin{eqnarray*}
\cosh(length(\beta))&=&\cosh(length(\beta_3))\sin(\pi/2-|\alpha|)\\
 &\ge & \cosh(2\delta)\cos(\delta) > \cosh(C). 
\end{eqnarray*}
So $length(\beta) > C$. Since $length(\beta) \le length(\beta_3) \le 4\delta$ this proves the claim.

For $i=1,2$ choose nonintersecting lifts $\tgamma_i \subset T_1(\H^2)$ of $\gamma_i$
so that the perpendicular segment $\beta$ from $\tgamma_1$ to
$\tgamma_2$ is on the left sides of both $\tgamma_1$ and
$\tgamma_2$. By the claim we may assume that $C \le length(\beta) \le 4\delta$. Let
\begin{eqnarray*}
M_j =2e^{-length(\gamma_j)/4}.
\end{eqnarray*}
Since $length(\gamma_j)-L\in (-\epsilon,\epsilon)$ it must be that $M_1+M_2 \le C$. So there exists a geodesic $\sigma \subset \H^2$ between $\tgamma_1$ and
$\tgamma_2$ so that for $i=1,2$, the shortest path from $\sigma$ to
$\tgamma_j$ has length
$M_j$. See figure \ref{fig:K}.

\begin{figure}[htb] 
\begin{center}
 \ \psfig{file=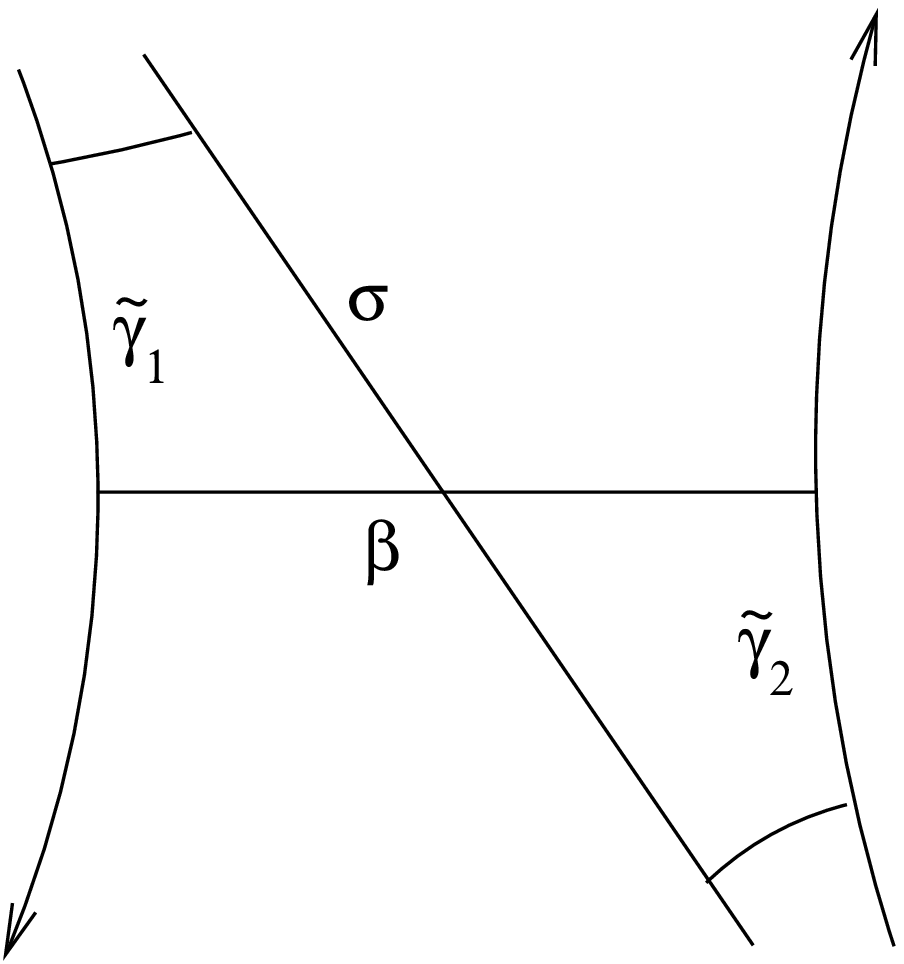,height=2.2in,width=2.2in}
 \caption{$\tgamma_1,\tgamma_2, \beta, \sigma$.}
 \label{fig:K}
 \end{center}
 \end{figure}

Let $\sK$ be the skew right-angled hexagon with three alternating sides contained in $\tgamma_1$, $\sigma$ and $\tgamma_2$. Let $s$ be the length of the subsegment
of $\sigma$ contained in this hexagon. By the law of cosines,
\begin{eqnarray*}
\cosh(s)&=&\frac{\cosh(length(\beta))-\cosh(M_1)\cosh(M_2)}{\sinh(M_1)\sinh(M_2)}\\
        &=&e^{length(\gamma_1)/4 + length(\gamma_2)/4}(\cosh(length(\beta))-1)/4 +O(1).
\end{eqnarray*}
So
\begin{eqnarray*}
s = L/2 + (\rho_1+\rho_2)/4 + \log(\cosh(length(\beta))-1) -\log(2) + O(e^{-L/2})
\end{eqnarray*}
where $\rho_i = length(\gamma_i)-L$ for $i=1,2$. 

Using a variant of the construction from section
\ref{section:isometry} we will obtain an isometry $h \in \Gamma$ whose
axis is very close to $\sigma$ such that $tr.length(h) \in
(L-\epsilon,L+\epsilon)$. Let $w$ be the unit vector tangent to $\sigma$ so that the basepoint of $w$ is a distance $s/2$ from the shortest path between $\sigma$ and either $\tgamma_1$ or $\tgamma_2$.

%half-way in between $\tgamma_i \cap \alpha_i$ where $\alpha_i$ 

%For simplicity, assume
%that $\sigma$ is the geodesic with endpoints $\{0,\infty\}$ in the
%upper half plane model and, for $j=1,2$, the point $p_j=e^{L/2 +
%  (2j-3)s/2}i \in \H^2$ is contained in the common
%perpendicular between $\sigma$ and $\gamma_j$. Let $a=(i,i), b=(e^{L}i,e^{L}i) \in T_1(\H^2)$.

 The right action of $\{V_t: t\in\R\}$ on $\Gamma \backslash PSL_2(\R)$
 is a horocycle flow. It well known \cite{Hed} that every orbit of
 this flow is dense in $\Gamma \backslash PSL_2(\R)$. Since $S$ is
 compact, there exists a time $T_0>0$ (depending only on $S$) such
 that for any vectors $x,y\in T_1(S)$ there
 exists a time $T$ with $|T|<T_0$ such that $y\in x V_T \sB_{\epsilon/5}$ (where $\sB_\epsilon$ is defined in section
 \ref{section:notation}). In particular, there exists a time $T$ with
 $|T| < T_0$ and a $h \in \Gamma$ such that  $hwG_{-L/2} = wG_{L/2}V_T Mat[a,b,c,d]$ for some $a,b,c,d\in \R$ with $|\ln(a)|,|b|,|c|,|d-1| < \epsilon/5$.  Note

\begin{displaymath}
V_T Mat[a,b,c,d]:=
\left[ \begin{array}{cc}
a & b \\
Ta +c & Td +b
\end{array} \right].
\end{displaymath}
By proposition \ref{pro:translation lengths} the translation length of $h$ (which is conjugate to $G_L V_T Mat[a,b,c,d]$) equals $L + 2\ln(a) + O(e^{-L})$. Since $\ln(a)<\epsilon/5$, if $L$ is
sufficiently large, the translation length of $h$ will be in $(L-\epsilon/2,L+\epsilon/2)$.

For $j=1,2$, let $g_j \in \Gamma$ be the isometry with axis $\tgamma_j$ and translation length equal to the length of $\gamma_j$. We will show that for $j=1,2$, $\{g_j,h\}$ determines a pair of pants $P_j \in \P_L(1,1,1)$ so that $(P_1,\pi(axis(h)), P_2)$ determines an $L$-clothesline of length 2 from $\gamma_1$ to $\gamma_2$. 

Lemma \ref{lem:x} above implies that for $i=1,2$ the length of the shortest path from $axis(h)$ to $\tgamma_i$ has length
\begin{eqnarray*}
M'_i := M_i +  O\Big(\frac{|b|+|Ta + c|+1}{a}\sqrt{\cosh(length(\beta))-1)}e^{-L/4 + (\rho_1+\rho_2)/8}\Big).
\end{eqnarray*}
Thus if $length(\beta)$ is sufficiently
small and $L$ is sufficiently large then 
\begin{eqnarray*}
2e^{-length(\gamma_i)/4 - \epsilon/8} \le M'_i \le 2e^{-length(\gamma_i)/4 + \epsilon/8} .
\end{eqnarray*}
By lemma \ref{lem:pants estimates}, this implies that $\{g_i,h\}$ determines a pair of pants $P_i \in \P_L(1,1,1)$ as
claimed. Since $length(\beta) \le 4 \delta$ this implies the lemma.

\end{proof}

\begin{lem}
Assume $Comm(\Gamma) < Isom^+(\H^2)$. Let $\delta > 0$. If $L$ is sufficiently large, then for every $\gamma_1 \in \G_L$, there exists an $L$-clothesline of length 1 between $\gamma_1$ and a $\delta$-dense geodesic $\gamma_2 \in \G_L$.
\end{lem}

\begin{proof}
For $T>0$, we say that a vector $w \in T_1(S)$ is $(T,\delta)$-dense if the geodesic segment $\{wG_t: 0 <  t < T\}$ is $\delta$-dense.

By ergodicity of the geodesic flow, for some $T>0$, there exists a $(T,\delta/2)$-dense vector $v_1 \in T_1(S)$. There exists a constant $C=C(T,\delta/2)$ such that any vector $w$ that has distance at most $C$ from $v_1$ is $(T,\delta)$-dense.

%Let $v_2= v_1 G_T$.

Let $L>>0$ and let $\gamma_1 \in \G_L$. Let $\tgamma_1$ be a lift of $\gamma_1$ to $\H^2$. Let $M=2e^{-length(\gamma_1)/4}$. Let $w \in T_1(\H^2)$ be a unit vector with basepoint of distance $M$ from $\tgamma_1$. Assume that $w$ is on the left side of $\tgamma_1$ and the shortest path from the basepoint of $w$ to $\tgamma_1$ is perpendicular to $w$. For $t\in \R$ let $w_t = wR_{-\pi/2}G_{t}R_{\pi/2}$.

 By the equidistribution theorem \ref{thm:equidistribution}, there exists a $T_0 > 0$ (depending only on the surface $S$ and $C$) such that for all $L$ sufficiently large, there exist times $t_1,t_2 $ such that 

\begin{itemize}
\item $0 \le t_1 \le t_2 \le T_0 e^{-L/3}$
\item $\pi\big(w_{t_1} G_{L/3}\big)$ is within a distance $C/2$ of $v_1$
\item $\pi\big(w_{t_2}G_{-L/2}\big) \in \pi\big(w_{t_2}G_{L/2}\big) \sB_{\epsilon/2}$.
\end{itemize} 
Here $\pi:T_1(\H^2) \to T_1(S)$ is the quotient map. Let $h = g_L(w_{t_2})$. A standard calculation shows that there is a vector $u'$ in $axis(w_{t_2})$ that is within a distance $O(e^{-L/6})$ of $w_{t_1}G_{L/3}$. It follows from lemma \ref{lem:x} that there is a vector $u \in axis(h)$ that is within a distance $O(e^{-L/3})$ from $u'$. Thus $\pi(u)$ has distance $C/2 + O(e^{-L/6})$ from $v_1$. If $L$ is large enough, this is less than $C$. So we may assume $\pi(axis(h))$ is $\delta$-dense. 

From the definition of $h$ and proposition \ref{pro:translation lengths} it follows that $tr.length(h) \in (L-\epsilon,L+\epsilon)$. From lemma \ref{lem:x}, it follows that the shortest path between $axis(h)$ and $\tgamma_1$ has length $M+O(e^{-L/3})$. Let $g \in \Gamma$ have axis equal to $\tgamma_1$ and translation length equal to the length of $\gamma_1$. It now follows from lemma \ref{lem:pants estimates} that if $L$ is large enough then the immersion $j:P \to S$ determined by $\{g,h\}$ is in $\P_L(1,1,1)$. To finish the lemma, set $\gamma_2 = \pi(axis(h))$.

\end{proof}

\begin{remark} The restriction that $Comm(\Gamma)<Isom^+(\H^2)$ can be removed by employing the variant of the isometry construction given in the previous lemma instead. 
\end{remark}

\begin{lem}
Given any $\delta >0$
\begin{eqnarray*}
\lim_{L \to \infty} \frac{\#\{\gamma \in \G_L : \gamma \textrm{ is $\delta$ dense}\}}{|\G_L|} = 1.
\end{eqnarray*}
\end{lem}

\begin{proof}
The set of $(T,\delta)$-dense vectors in $T_1(S)$ is open for every $T$ and $\delta$. By ergodicity of the geodesic flow,
\begin{eqnarray*}
\lim_{T \to \infty} \lambda( w \in T_1(S): w \textnormal{ is $(T,\delta)$-dense}) =1
\end{eqnarray*}
where $\lambda$ is Haar probability measure on $T_1(S)$. The lemma now follows from theorem \ref{thm:bowen}.

\end{proof}

Theorem \ref{thm:clotheslines} follows immediately from the preceding three lemmas.

\section{Calculations}\label{section:calculations}

\subsection{Translation Length}

If $g \in Isom^+(\H^2)$, then the translation length of $g$ is the shortest possible distance between a point $p \in \H^2$ and its translate $gp$ (minimized over all $p \in \H^2$). More generally, for $g \in Isom^+(\H^3)$, we define the displacement $\mu(g)$ by $\cosh(\mu(g)/2)=trace(g)/2$ and $Re(\mu(g))\ge 0$ where $g$ is identified with its matrix representative in $PSL_2(\C)$. It can be shown that $\mu(g)=tr.length(g)$ if $g \in Isom^+(\H^2)$. The following is a short calculation.
\begin{pro}\label{pro:translation lengths}
If $g=G_L Mat[a,b,c,d]$ then
\begin{equation}\label{eqn:g}
g=\left[ \begin{array}{cc}
a e^{L/2} &be^{L/2}  \\
ce^{-L/2}& d e^{-L/2} 
\end{array} \right]
\end{equation}
and the displacement of $g=G_L Mat[a,b,c,d]$ satisfies
\begin{displaymath}
\mu(g) =  L + 2\ln(a) + O(e^{-L}).
\end{displaymath}
\end{pro}
If $g=g_L(w)$ is as defined as in section \ref{section:isometry} so that $gwG_{-L/2}=wG_{L/2}Mat[a,b,c,d]$ then $g$ is conjugate to $G_L Mat[a,b,c,d]$. So its translation length is given by the above.

\subsection{Calculations for hexagon $\sH$}\label{section:hexagon H}

\begin{figure}[h]
 \begin{center}
 \ \psfig{file=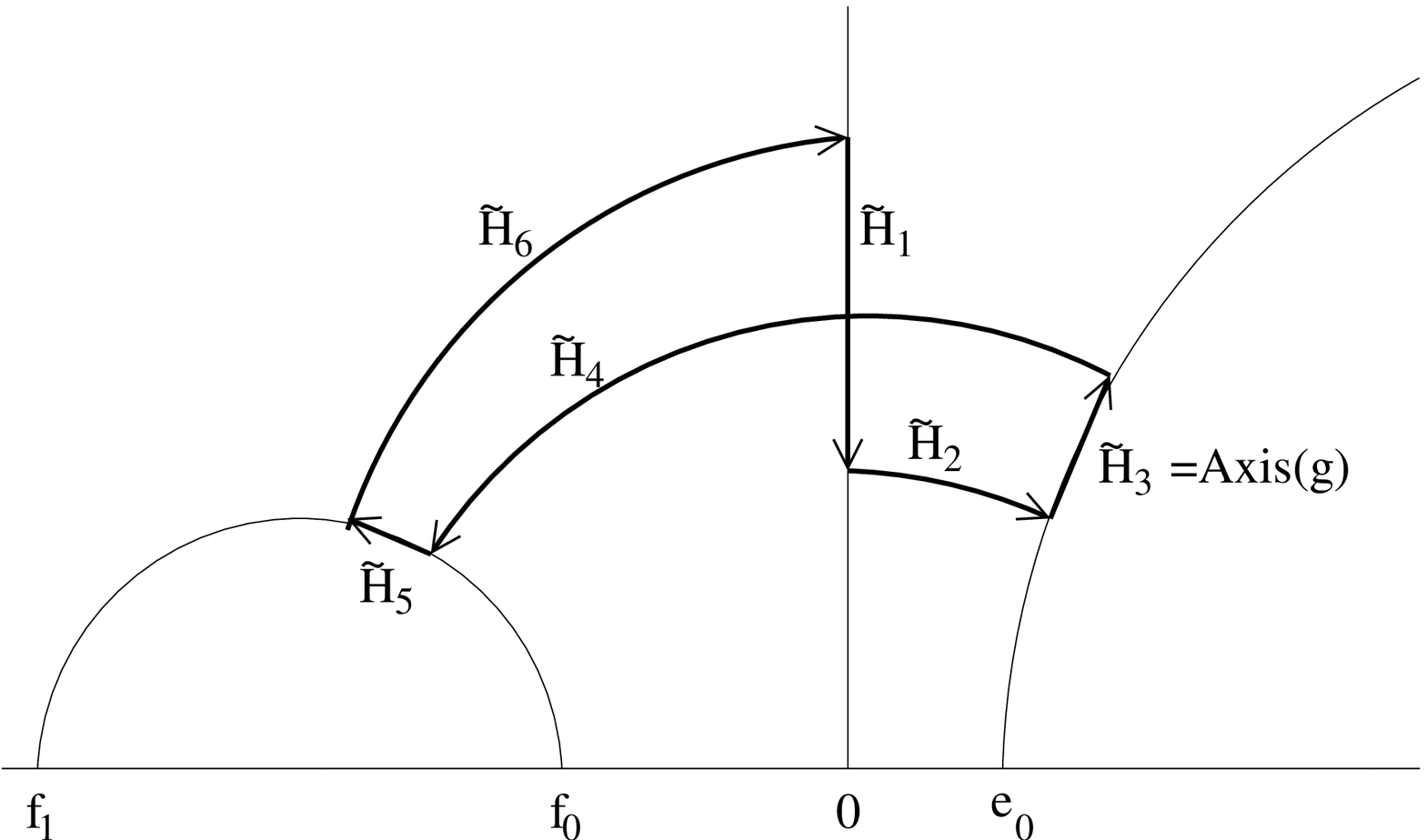, height=3in, width=4.75in}
 \caption{The hexagon $\hH$ (in the upperhalf space model).}
 \label{fig:hexagonH}
 \end{center}
 \end{figure}

The goal of this section is to prove estimates on the sidelengths (and widths) of a right-angled hexagon $\sH$ related to the construction of isometries in section \ref{section:isometry}. From these estimates corollary \ref{cor:2d}, corollary \ref{cor:perp} and lemma \ref{lem:x} will follow. By a right-angled hexagon $\sH$ we mean an ordered  $6$-tuplet of oriented geodesics $(\tH_1,\tH_2,..,\tH_6)$ of $\H^3$ such that $\tH_i$ is orthogonal to $\tH_{i+1}$ and $\tH_i$ is not equal to $\tH_{i+2}$ for any $i$ modulo $6$. For the sake of generality, we allow $\sH$ to be nonplanar. To handle this we need an appropriate substitute for sidelength called width (see subsection \ref{subsection:doublecross}). Let $H_i=\mu(\tH_{i-1},\tH_{i+1};\tH_i)$ be the width of the double-cross $(\tH_{i-1},\tH_{i+1};\tH_i)$. To define $\sH$, let
\begin{eqnarray*}
g = G_LMat[a,b,c,d]
\end{eqnarray*}
where $L>0$, $a,b,c,d\in \C$ such that $ad-bc=1$ and $Re(a)>0$. Let $J,M:(0,\infty)\to \C$ be functions of the parameter $L$. In most applications, $J=e^{L/2}$ and $M=m_1e^{-m_2L}$ for some $m_1, m_2 >0$. $\sH=(\tH_1,...,\tH_6)$ is determined by:
\begin{enumerate}
\item $\tH_1$ is equal to the geodesic with endpoints $\{0,\infty\}$
  oriented from $\infty$ to $0$.
\item $\tH_2$ is oriented from $\tH_1$ to $\tH_3$ if they do not
  intersect. Otherwise $\tH_2$ is oriented so that if $v_i$ is a unit
  vector based at the intersection point $x=\tH_1 \cap \tH_3$ pointing
  in the direction of $\tH_i$ ($i=1,2,3$) then $(v_1,v_3,v_2)$ is a
  positively oriented basis for the tangent space at $x$ in
  $\H^3$. 
\item $\tH_3$ is equal to $axis(g)$ oriented from $e_0$ to $e_1$ (defined below).
\item $\tH_4$ is oriented from $\tH_3$ to $\tH_5$ if they do not
  intersect. Otherwise $\tH_4$ is oriented so that if $v_i$ is a unit
  vector based at the intersection point $x=\tH_3 \cap \tH_5$ pointing
  in the direction of $\tH_i$ ($i=3,4,5$) then $(v_3,v_5,v_4)$ is a
  positively oriented basis for the tangent space to at $x$ in
  $\H^3$. 
%\item $\tH_5$ is oriented arbitrarily.
\item $\tH_6$ has endpoints $\pm J$ and is oriented from $-J$ to $+J$.
\item $H_6=\gM + i\pi$.
%\item $\tH_i$ is oriented from its intersection point with $\tH_{i-1}$ to its intersection point with $\tH_{i+1}$ assuming these do not coincide, which will be the case in all applications.
\end{enumerate}
Hexagon $\sH$ is depicted in figure \ref{fig:hexagonH}. 
%set the M versus \gM notation straight

Let $j=|L/2 - \log(|J|)|$. So $j$ is the hyperbolic distance from the point $e^{L/2}i\in \H^3$ to the point where $\tH_6$ intersects $\tH_1$. Let 
\begin{eqnarray*}
E=(\frac{|c|+|b|+|bc| + 1}{a})e^{j-L/2}).
\end{eqnarray*}
 The main theorem is:
\begin{thm}\label{thm:main}
Assume $j < L/2$ and $|M| > |E|$. Then there exists $(\sigma_1,\sigma_2) \in \{(1,1),(-1,0)\}$ such that
\begin{eqnarray*}
\wH_4 &=&  \sigma_1 M + \sigma_1 x + O(\coth(M)E^2),\\
\wH_5 &=&  \sigma_2 i\pi + \frac{y}{\sinh(M)} + O \Big(\frac{E^3}{\sinh^3(M)} \Big)
\end{eqnarray*}
where
\begin{eqnarray*}
x &=& \frac{c}{a} e^{-L}J + \frac{ -b}{a}(1/J) + O(Je^{-2L} + (1/J)e^{-L})\\
y^2 &=& \frac{4bc}{a^2} e^{-L} + x^2 + O(e^{-2L}).
\end{eqnarray*}
\end{thm}
This theorem is proven in subsection \ref{section:main}. The next lemma determines $H_4$ and $H_5$ more fully in special cases.
\begin{lem}\label{lem:signs}
Assume $a,b,c,d \in \R$. If $L$ is sufficiently large, the following hold.
\begin{itemize}
\item If $M > |E|$ and $J=e^{L/2}$ then
\begin{eqnarray*}
\wH_4&=& M+ \frac{c-b}{a} e^{-L/2} + O(e^{-3L/2}+ \coth(M)E^2)\\
\wH_5&=& \frac{c+b}{a} e^{-L/2} + O\Big(e^{-3L/2}+ \frac{E^3}{\sinh^3(M)} \Big).
\end{eqnarray*}
\item If $M=-i\pi/2$ and $J=-ie^{L/2}$ then
\begin{eqnarray*}
\wH_4 &=& i\pi/2 +i \frac{c+b}{a} e^{-L/2} + O(e^{-3L/2}+ \coth(M)E^2)\\
\wH_5 &=&\frac{c-b}{a} e^{-L/2} + O\big(e^{-3L/2}+ \frac{E^3}{\sinh^3(M)} \big).
\end{eqnarray*}
\end{itemize}
\end{lem}
We defer the proof of this lemma to the end of section \ref{section:main}. We can now prove corollary \ref{cor:2d} and corollary \ref{cor:perp}.
\begin{proof}(of corollary \ref{cor:2d})
 After conjugating we may assume $w=(e^{L/2}i,e^{L/2}i)$ in the upperhalf space model. Let $F_1,F_2,F_3$ be as given after the statement of corollary \ref{cor:2d}. Apply theorem \ref{thm:main} with $J=e^{L/2}$ and $M=M$. So $\tH_1$ is $axis(w)$ with reverse orientation, $\tH_3=axis(g)$ and $\tH_5=\tgamma$. $H_4$ is now the distance between $axis(w)$ and $\gamma$. $H_5$ is the signed distance from $p$ to $q$. The corollary now follows from the above lemma and proposition \ref{pro:translation lengths}.
\end{proof}
\begin{proof}(of corollary \ref{cor:perp})
  After conjugating we may assume $w=(e^{L/2}i,e^{L/2}i)$ in the upperhalf space model. Set $M$ in theorem \ref{thm:main} equal to $-i\pi/2$ and $J=-ie^{L/2}$. So $\tH_1=axis(w)$ with reversed orientation, $\tH_3=axis(g)$ and $\tH_5=axis(w^\perp)$. $H_4/i$ is the angle between $axis(w^\perp)$ and $axis(g)$. $H_5$ is the signed distance between the basepoint of $w$ and $axis(w^\perp)\cap axis(g)$. The corollary now follows from the above lemma and proposition \ref{pro:translation lengths}.
\end{proof}

Before proving the main theorem, we prove some identities, as opposed to estimates, regarding the widths of $\sH$. For this, we need a formula for the fixed points of $g$. In general the fixed points of
\begin{displaymath}
\left[ \begin{array}{cc}
a & b \\
c & d
\end{array} \right] \in PSL_2(\C)
\end{displaymath}
acting by fractional linear transformations on $\C$ are 
\begin{eqnarray}\label{eqn:fixed points}
\{ \frac{ a - d \pm \sqrt{(a+d)^2 - 4} }{ 2c} \}&& \textnormal{if $c \ne 0$}\\
\{\infty, \frac{b}{d-a}\} &&  \textnormal{if $c = 0$}.
\end{eqnarray}
It follows from equation \ref{eqn:g} (subsection
\ref{section:isometry}) that (if $c \ne 0$), $e_0=(N_1-N_2)/D$ and $e_1=(N_1+N_2)/D$ are the fixed points of $g$ where 

\begin{eqnarray*}
N_1 &=& ae^{L/2} -de^{-L/2}.\\
N_2^2 &=& (ae^{L/2} +de^{-L/2})^2 - 4.\\
D   &=&  2ce^{-L/2}.
\end{eqnarray*}
We choose $N_2$ to be a root of the above equation with nonnegative real part. When $L$ is large and $Re(e^a)>0$, $e_0$ is close to zero and $e_1$ is ``close'' to $\infty$.

\begin{thm}\label{thm:identities}
\begin{eqnarray*}
\cosh(\wH_2) &=& -N_1/N_2\\
\cosh(\wH_4) &=& (N_1/N_2)\cosh(\gM)+ x\sinh(\gM)\\
\wH_6 &=& M + i\pi.
\end{eqnarray*}

where 
\begin{eqnarray*}
x &=& \frac{1}{2}\frac{DJ^2 + (N_1^2 - N_2^2)/D }{J N_2 }.
\end{eqnarray*}
\end{thm}

\begin{remark} No hypotheses on $M$ or $J$ are needed above. Assuming the theorem, we can use the law of cosines and the law of sines to determine formulas for the other three widths of $\sH$.
\end{remark}
We prove this theorem first and then derive theorem \ref{thm:main} using it.

\subsubsection{Trigonometry Background: The Cross Ratio}\label{subsection:cross ratio}

The cross ratio $R$ of $(a,b,c,d) \in \C^4$ is defined by

\begin{equation}\label{eqn:cross ratio}
R(a,b,c,d) = \frac{ (a-c)(b-d) }{(a-d)(b-c) }.
\end{equation}

The cross ratio is invariant under the action of $PSL_2(\C)$ by fractional linear transformations on $\C$. Note that
\begin{displaymath}
R(b,a,c,d) = \frac{1}{R(a,b,c,d)}
\end{displaymath}
and $R(a,b,c,d) = R(c,d,a,b)$. $R$ is extended to $(\C\cup\{\infty\})^4$ in the natural way. For example, $R(\infty, b,c,d) = \frac{b-d}{b-c}$. 

\subsubsection{Trigonometry Background: Double Crosses}\label{subsection:doublecross}
The material in this subsection is detailed more thoroughly in \cite{Fen}. Suppose $u, u'$ are the endpoints of a geodesic $\gamma_1$ oriented from $u$ to $u'$ and $v,v'$  are the endpoints of a geodesic $\gamma_2$ oriented from $v$ to $v'$. Suppose also that $\gamma_3$ is a geodesic perpendicular to both $\gamma_1$ and $\gamma_2$. Let $w,w'$ be the endpoints of $\gamma_3$ which we assume is oriented from $w$ to $w'$. The triple $(\gamma_1,\gamma_2; \gamma_3)$ is called a {\bf double cross}. We define the {\bf width} $\mu(\gamma_1,\gamma_2;\gamma_3) = \mu \in \C/<2\pi i>$ of the double cross by the equation
\begin{eqnarray}\label{eqn:R1}
exp(\mu) &=& R(u,v,w',w)=-R(u',v,w',w) \\
         &=& -R(u,v',w',w)\\
         &=& R(u',v',w',w).
\end{eqnarray}
So
\begin{displaymath}\begin{array}{ll}
\exp(\mu(\gamma_1,\gamma_2;\gamma_3)) &= R(u,v,w',w)\\
                                &= 1/R(v,u,w',w)\\
                                &=\exp(-\mu(\gamma_2,\gamma_1;\gamma_3 )).
\end{array}
\end{displaymath}
Hence $\mu(\gamma_1,\gamma_2;\gamma_3) =
-\mu(\gamma_2,\gamma_1;\gamma_3)$. $\mu$ also satisfies the equation

\begin{equation}\label{eqn: displacement}
R(u,u',v,v') = \tanh^2(\mu/2).
\end{equation}

The latter equation determines $\mu$ only up to a sign. If we denote $\gamma_i$ with the opposite orientation by $-\gamma_i$, then we have

\begin{displaymath}
\mu(\gamma_1,\gamma_2;-\gamma_3) = -\mu(\gamma_1,\gamma_2;\gamma_3)
\end{displaymath}

and

\begin{displaymath}
\mu(-\gamma_1,\gamma_2;\gamma_3) = \mu(\gamma_1,\gamma_2;\gamma_3) + i\pi.
\end{displaymath}

 The real part of $\mu$ is the signed distance between $\gamma_1$ and $\gamma_2$. The imaginary part measures the amount of turning between $\gamma_1$ and $\gamma_2$. To be precise, $\mu$ is the displacement of the isometry $g$ with oriented axis $Axis(g)=(w,w')$ such that $g \gamma_1=\gamma_2$. Suppose as above that
\begin{displaymath}\label{eqn:tanh}
R=R(u,u',v,v') = \tanh^2(\mu/2).
\end{displaymath}
Then
\begin{eqnarray*}
\frac{1+R}{1-R} &=&\cosh(\mu)\\
\frac{\sqrt{R}}{1-R} &=&(1/2)\sinh(\mu).
\end{eqnarray*}

%\subsubsection{Law of Cosines and Sines}

%If $(S_1,S_2,S_3,S_4,S_5,S_6)$ is an {\it ordered} $6$-tuplet of oriented geodesics such that $S_i$ is orthogonal to $S_{i+1}$ and $S_i$ is not equal to $S_{i+2}$ for any $i$ modulo $6$, then it is called a {\bf right angled hexagon}. By orthogonal, we will mean that $S_i$ and $S_{i+1}$ intersect in $\H^3$ at a right-angle (this constrasts a little with \cite{Fen} where the word ``normal'' is used to allow the possibility that $S_i$ and $S_{i+1}$ share an endpoint at infinity). In the terminology of \cite{Fen} all the side-lines that we allow are proper.

%Similarly, if $(S_1,S_2,S_3,S_4,S_5)$ is an {\it ordered} $5$-tuplet of oriented geodesics such that $S_i$ is orthogonal to $S_{i+1}$ and $S_i$ is not equal to $S_{i+2}$ for any $i$ modulo $5$, then it is called a {\bf right angled pentagon}. 

The following lemma is classical. It appears in \cite{Fen}.

\begin{lem}\label{lem:right angled hexagons}
Let $\sH=(\tH_1,..., \tH_6)$ be a right-angled hexagon. Let $H_i = \mu(\tH_{i-1}, \tH_{i+1}; \tH_i)$ denote the width of the double cross $(\tH_{i-1},\tH_{i+1};\tH_i)$. Then the following relations hold.

\begin{enumerate}
\item The law of sines: $\frac{\sinh(H_1)}{\sinh(H_4)} = \frac{\sinh(H_3)}{\sinh(H_6)} = \frac{\sinh(H_5)}{\sinh(H_2)}$.

\item The law of cosines: $\cosh(H_i) = \cosh(H_{i-2}) \cosh(H_{i+2}) + \sinh( H_{i-2}) \sinh( H_{i+2} ) \cosh(H_{i+3})$ for all $i$ with indices considered modulo $6$.
\end{enumerate}
\end{lem}

%fix the props below to remove the estimates.
\subsubsection{Proof of theorem \ref{thm:identities}}\label{section:identities}

\begin{pro}\label{pro:H2}
If $c\ne 0$
\begin{displaymath}\begin{array}{ll}
\cosh(\wH_2) &= -N_1/N_2.
\end{array}
\end{displaymath} 

\end{pro}

\begin{proof}
Recall the definition of the cross ratio $R$ (subsection \ref{subsection:cross ratio}). Let
\begin{displaymath}
R:=R(\infty,0,e_0,e_1) = \frac{e_1}{e_0} = \frac{N_1+N_2}{N_1-N_2}.
\end{displaymath}
Recall that $\wH_2 = \mu(\tH_1,\tH_3;\tH_2)$. By subsection \ref{subsection:doublecross}, $\tanh^2(\mu(\tH_1,\tH_3;\tH_2)/2)=R$. So
\begin{eqnarray*}
\cosh(\wH_2) = \frac{1+R}{1-R} = \frac{-N_1}{N_2}.
\end{eqnarray*}
\end{proof}

\begin{pro}
Let $\{f_0,f_1\}$ be the endpoints of $\tH_5$ such that $\tH_5$ is
oriented from $f_0$ to $f_1$. Let $\m=\tanh(\gM/2)$. Then
\begin{displaymath}\begin{array}{ll}
f_0 &= -\m J\\
f_1 &= -J/\m.
\end{array}
\end{displaymath}

\end{pro}

\begin{proof}
 By
equation \ref{eqn:R1},
\begin{eqnarray*}
e^{H_6} &=& R(f_0,\infty,J,-J)\\
        &=& \frac{f_0-J}{f_0+J}.
\end{eqnarray*}
But $e^{H_6} = e^{M + i\pi} = -e^M$. Solving for $f_0$ yields $f_0=
-J\tanh(M/2)=-\m J$. By equation \ref{eqn:R1},
\begin{eqnarray*}
e^{H_6} &=& -R(f_1,\infty,J,-J)\\
        &=& -\frac{f_1-J}{f_1+J}.
\end{eqnarray*}
Solving for $f_1$ yields $f_1=-J/\m$.

\end{proof}

\begin{pro}\label{pro:H4}
If $c\ne 0$
\begin{eqnarray*}
\cosh(\wH_4) &=& (N_1/N_2)\cosh(\gM)+ \sinh(\gM)x;
\end{eqnarray*}
where 
\begin{eqnarray*}
x &=& (1/2) \frac{DJ^2 + (N_1^2 - N_2^2)/D }{J N_2 }.
\end{eqnarray*}

\end{pro}

\begin{proof}

We compute $R(f_0,f_1,e_0,e_1)$ as follows.

\begin{eqnarray*}
R(f_0,f_1,e_0,e_1) &=& \frac{(-\m J -\frac{N_1-N_2}{D}) (-J/\m -\frac{N_1+N_2}{D}) }{ (-\m J -\frac{N_1+N_2}{D}) (-J/\m -\frac{N_1-N_2}{D}) }\\
                   &=&\frac{(-D\m J - (N_1 - N_2))(-DJ - \m(N_1 + N_2)) }{ (-D \m J - (N_1 +N_2)  )( -DJ - \m(N_1 - N_2) ) }\\
                   &=&\frac{ D^2\m J^2 + DJ(N_1 - N_2) + D \m^2 J(N_1 + N_2) + \m(N_1^2 - N_2^2) }{ D^2\m J^2 + DJ(N_1 + N_2) + D \m^2 J(N_1 - N_2) + \m(N_1^2 - N_2^2)}\\
                   &=&\frac{ D\m J^2 + J(N_1 - N_2) +  \m^2 J(N_1 + N_2) + \m(N_1^2 - N_2^2)/D }{ D\m J^2 + J(N_1 + N_2) +  \m^2 J(N_1 - N_2) + \m(N_1^2 - N_2^2)/D} = (A/B)
 %&\\
                   %&=\m^2 + \frac{ D(\m-\m^3) J^2 + (1-\m^4) J(N_1 - N_2) + (\m-\m^3)(N_1^2 - N_2^2)/D }{ D\m J^2 + J(N_1 + N_2) +  \m^2 J(N_1 - N_2) + \m(N_1^2 - N_2^2)/D}\\
\end{eqnarray*}
where $A$ is the numerator in the line above and $B$ is the denominator. Since $R = \tanh^2(\wH_4/2)$ we obtain the following.

\begin{eqnarray*}
\cosh(\wH_4) &=&\frac{1 + R}{1-R} = \frac{1 + (A/B)}{1- (A/B)} = \frac{B+A}{B-A}\\
%             &=&\frac{ 2D\m J^2 + (1+\m^2)J(N_1 - N_2) +  (1+\m^2) J(N_1 + N_2) + 2\m(N_1^2 - N_2^2)/D}{ (\m^2-1)J(N_1 - N_2) +  (1-\m^2) J(N_1 + N_2)}\\
             &=&\frac{2D\m J^2 + 2(1+\m^2)J N_1 + 2\m(N_1^2 - N_2^2)/D}{ 2(1-\m^2)J N_2}\\
      %       &=&\frac{D\m J^2 + (1+\m^2)J N_1 + \m(N_1^2 - N_2^2)/D}{ (1-\m^2)J N_2}\\
             &=&\frac{(1+\m^2)J N_1}{ (1-\m^2)J N_2} + \frac{D\m J^2  + \m(N_1^2 - N_2^2)/D}{ (1-\m^2)J N_2}\\
             &=&\cosh(\gM)(N_1/N_2) + \frac{\m}{1-\m^2} \left(\frac{DJ^2 + (N_1^2 - N_2^2)/D}{J N_2}\right)\\
             &=&\cosh(\gM)(N_1/N_2) + (1/2)\sinh(\gM)\left(\frac{DJ^2 + (N_1^2 - N_2^2)/D}{J N_2}\right)\\
&=& \cosh(\gM)(N_1/N_2) + x\sinh(\gM).
\end{eqnarray*}
\end{proof}

If $c\ne 0$, theorem \ref{thm:identities} follows from the three propositions above. The case $c=0$ follows from taking a limit of the $c\ne 0$ case.

\subsection{Proof of theorem \ref{thm:main}}\label{section:main}

\begin{pro}\label{pro:fixed points}
%If $|T|,|\nu|,|\delta|,|\theta|<\epsilon$ then
%\begin{eqnarray*}
%\cosh(H_2) &=& -N_1/N_2 = -1 + O(\epsilon e^{-L}).\\
%x &=& O(\epsilon j_1 e^{(j_2-1)L} + \frac{\epsilon}{j_1} e^{-j_2L}).
%%y &:=& \frac{D^2 J^4 + (N_1^2 - N_2^2)^2/D^2 }{J^2 N_2^2} = O(e^{2(x-1)L} + e^{-2xL}).
%\end{eqnarray*}

\begin{eqnarray*}
\cosh(H_2) &=& -N_1/N_2 = -1 + \frac{2bc}{a^2} e^{-L} + O(e^{-2L}).\\
x &=& \frac{c}{a} e^{-L}J + \frac{ -b}{a}(1/J) + O(Je^{-2L} + (1/J)e^{-L}).
%y &=& 
\end{eqnarray*}
In particular, $x= O(E)$ and $\cosh(H_2) = -1 + O(E^2)$ if $j<L/2$.

%Here the constant implicit in the $O(\cdot)$ notation is independent of all constants in the construction ($L,T,\nu,\delta,\theta, \epsilon,j_1,j_2,m_1,m_2$).
\end{pro}

\begin{proof}
Recall:
\begin{eqnarray*}
N_1 &=& ae^{L/2} -d e^{-L/2}.\\
N_2^2 &=& (ae^{L/2} +d e^{-L/2})^2 - 4.\\
D   &=&  2ce^{-L/2}.
\end{eqnarray*}
For $z$ close to zero, $\sqrt{1-z}=1 -(1/2)z + O(z^2)$. So if $z$ is very large, $\sqrt{z^2 -4} = z\sqrt{1 - 4/z^2} = z - 2/z + O(1/z^3)$. Hence,

\begin{eqnarray*}
N_2 &=& ae^{L/2} + d e^{-L/2} - \frac{2}{ae^{L/2} + d e^{-L/2}} + O(e^{-3L/2})\\
    &=& ae^{L/2} + (d-2/a)e^{-L/2} + O(e^{-3L/2}).
\end{eqnarray*}

So,
\begin{eqnarray*}
N_1-N_2 &=& -2(d-1/a)e^{-L/2} + O(e^{-3L/2})
\end{eqnarray*}
and
\begin{eqnarray*}
\frac{N_1-N_2}{N_2} &=& \frac{-2(d-1/a)e^{-L/2} + O(e^{-3L/2})}{ ae^{L/2} + (d-2/a)e^{-L/2} + O(e^{-3L/2})}\\
                    &=& \frac{-2(ad-1)}{a^2} e^{-L} + O(e^{-2L})\\
                    &=& \frac{-2bc}{a^2} e^{-L} + O(e^{-2L}).
\end{eqnarray*}
The estimate for $\cosh(H_2)$ now follows from proposition \ref{pro:H2}. Also, 
\begin{eqnarray*}
N_1^2-N_2^2 &=& (ae^{L/2}-de^{-L/2})^2 - (ae^{L/2}+de^{-L/2})^2 + 4\\
            &=& -4ad + 4=-4bc.
\end{eqnarray*}
Therefore,
\begin{eqnarray*}
\frac{N_1^2-N_2^2}{DN_2} &=& \frac{ -4bc}{2ce^{-L/2}( ae^{L/2} + (d-2/a)e^{-L/2} + O(e^{-3L/2}))  }\\
                      &=& \frac{ -2b}{a} + O(e^{-L}).
\end{eqnarray*}
Also,
\begin{eqnarray*}
D/N_2 &=& \frac{2ce^{-L/2}}{ae^{L/2} + (d-2/a)e^{-L/2} + O(e^{-3L/2})}\\
      &=& \frac{2c}{a} e^{-L} + O(e^{-2L}).
\end{eqnarray*}

Thus,
\begin{eqnarray*}
x &=& (1/2)DJ/N_2 + (1/2)(N_1^2-N_2^2)/(DJN_2)\\
 &=& \frac{c}{a} e^{-L}J + \frac{ -b}{a}(1/J) + O(Je^{-2L} + (1/J)e^{-L}).
\end{eqnarray*}

%When $J=e^{L/2}$, 
%\begin{eqnarray*}
%y^2 &=& 2(1-N_1/N_2) + x^2\\
%    &=&  \frac{4bc+4(c-b)^2}{a^2}  e^{-L} + O(e^{-2L})
%\end{eqnarray*}

\end{proof}

%**J sshould be a more intuitive quantity.

%Recall
%\begin{eqnarray*}
%E=(|c|+|b|+|bc|)e^{j-L/2}.
%\end{eqnarray*}
%Then $x=O(E)$, $E=O((|c|+|b|+|bc|)e^{j-L/2})$ and $N_1/N_2 -1 =O(E^2)$ if $j<L/2$.
\begin{pro}\label{pro:H4estimate}
Assume $|x| < |M|$ and $j<L/2$. Then there exists $\sigma_1 \in \{-1,+1\}$ such that:
\begin{eqnarray*}
\sinh(\wH_4) &=&  \sigma_1(N_1/N_2)\sinh(M) + \sigma_1 x\cosh(M) + O\Big( \frac{E^2}{\sinh(M)}\Big)\\
\wH_4 &=& \sigma_1 M + \sigma_1 x + O(\coth(M)E^2).
\end{eqnarray*}
\end{pro}

\begin{proof}
Recall that $\cosh(H_4)=(N_1/N_2)\cosh(M) + x \sinh(M)$. So,
\begin{eqnarray*}
\sinh^2(H_4) &=& \cosh^2(H_4)-1\\
              &=& \big( (N_1/N_2)\sinh(M) + x\cosh(M)\big)^2 - x^2 +(N_1/N_2)^2 - 1.
\end{eqnarray*}
Since $N_1/N_2 = 1 + O(E^2)$,
\begin{eqnarray*}
\sinh(H_4) &=& \sigma_1(N_1/N_2)\sinh(M) + \sigma_1 x\cosh(M) + O\Big( \frac{E^2}{\sinh(M)}\Big)
\end{eqnarray*}
for some $\sigma_1 \in \{-1,+1\}$. If $\sigma_1=1$ then 
\begin{eqnarray*}
\sinh(\wH_4 -\gM) &=& \cosh(M)\sinh(\wH_4) - \cosh(\wH_4)\sinh(\gM)\\
&=& (N_1/N_2)\cosh(\gM)\sinh(\gM) + \cosh^2(M)x\\
&& - (N_1/N_2)\cosh(M)\sinh(M) - \sinh^2(M)x + O(\coth(M)E^2 )\\
&=& x + O(\coth(M)E^2).
\end{eqnarray*}
In general, $\sinh(y)=y + O(y^3)$. So $H_4 = M + x + O(\coth(M)E^2)$. The calculation for $\sigma_1=-1$ is similar.
\end{proof}

\begin{pro}\label{pro:H5}
Assume $|x| < |M|$ and $j<L/2$. Then there exists $\sigma_2 \in \{1,0\}$ such that
\begin{eqnarray*}
\wH_5 &=& \sigma_2 i\pi  + \frac{y}{\sinh(M)} + O\Big( \frac{E^3}{\sinh^3(M)} \Big).
\end{eqnarray*}
where $y^2 = 1-(N_1/N_2)^2 + x^2$ and $(\sigma_1,\sigma_2) \in \{(1,1),(-1,0)\}$. $\sigma_1$ is as in the previous proposition.

\end{pro}

\begin{proof}
Recall $H_6=M+i\pi$. Assume $\sigma_1=+1$. By the law of cosines
\begin{eqnarray*}
\cosh(\wH_5) &=& \frac{\cosh(\wH_2) - \cosh(\wH_6)\cosh(\wH_4)} {\sinh(\wH_6)\sinh(\wH_4)}\\
             &=& \frac{-N_1/N_2 + \cosh(M)\big((N_1/N_2)\cosh(M) + x\sinh(M)\big)}{-\sinh(M)\sinh(H_4)}\\
             &=& -\frac{(N_1/N_2)\sinh(M) + x\cosh(M)}{\sinh(H_4)}\\
             &=& -1 + O(\cosh(M)E^2/\sinh^2(M) ).
\end{eqnarray*}
So $\wH_5 = i\pi + O(\sqrt{\cosh(M)}E/\sinh(M))$. Squaring, we obtain 
\begin{eqnarray*}
\cosh^2(\wH_5) &=& \frac{\big((N_1/N_2)\sinh(M) + x\cosh(M)\big)^2}{\big((N_1/N_2)\cosh(M) + x\sinh(M)\big)^2-1}\\
 &=& \frac{\big((N_1/N_2)\sinh(M) + x\cosh(M)\big)^2}{\big((N_1/N_2)\sinh(M) + x\cosh(M)\big)^2-1+(N_1/N_2)^2-x^2}\\
               &=& 1 + \frac{1-(N_1/N_2)^2 + x^2}{\sinh^2(H_4)}\\
&=& 1 + \frac{1- (N_1/N_2)^2+x^2}{\sinh^2(M)} +  O\Big( \frac{E^4}{\sinh^4(M)}\Big)\\
\end{eqnarray*}
So
\begin{eqnarray*}
\sinh(\wH_5) &=& \frac{y}{\sinh(M)} + O\Big( \frac{E^3}{\sinh^3(M)} \Big).
\end{eqnarray*}
where
\begin{eqnarray*}
y^2 = 1 - (N_1/N_2)^2 + x^2
\end{eqnarray*}
Since $\sinh(z)=z+O(z^3)$ the case $\sigma_1=1$ follows. The $\sigma_1=-1$ case is similar.

\end{proof}

\begin{proof}(of theorem \ref{thm:main})
The proof follows immediately from the preceding three propositions.
\end{proof}

\begin{proof}(of lemma \ref{lem:signs})
For now assume $M>|x|$ and $J=e^{L/2}$. Since $H_4$ is approximately $\sigma_1 M + \sigma_1 x$, $\tH_4$ has nonzero real part. So $\tH_3$ does not intersect $\tH_5$. Since $\tH_4$ is oriented from $\tH_3$ to $\tH_5$ it must have positive real part. Thus $\sigma_1=\sigma_2=1$. The equation for $H_4$ follows from the equation for $x$ in theorem \ref{thm:main}. Since $J=e^{L/2}$, a short calculation shows
\begin{eqnarray*}
y^2 = \frac{(c+b)^2}{a^2} e^{-L} + O(e^{-2L}).
\end{eqnarray*}
We test the above in the special case in which $(a,b,c,d)=(1,0,1,1)$. A short calculation shows that $e_0=0$ and $e_1=e^L - 1$. So figure \ref{fig:hexagonH} accurately depicts the incidences and orientations of the geodesics $\tH_1,..,\tH_6$ except that $\tH_1$ and $\tH_3$ are coincident and $\tH_2$ is degenerate. But this is a limit of cases like those depicted in the figure. The key observation is that $\tH_5$ is oriented from its intersection with $\tH_4$ to its intersection with $\tH_6$. Thus the real part of $H_5$ is positive. So in this case
\begin{eqnarray*}
y = \frac{(c+b)}{a} e^{-L/2} + O(e^{-3L/2}).
\end{eqnarray*}
By continuity, the above holds in general. This case now follows from the formula for $H_5$ in theorem \ref{thm:main}.

Now assume $M=-i\pi/2$ and $J=-ie^{L/2}$. A short calculation shows that
\begin{eqnarray*}
x&=&-i(c+b)e^{-L/2-a}\\
y^2&=&\frac{-(c-b)^2}{a^2} e^{-L}.
\end{eqnarray*}
In particular $x$ and $y$ are both purely imaginary. Since $\tH_5$ is orthogonal to both $\tH_1$ and $\tH_6$ it must have endpoints equal to $\pm e^{L/2}$. The orientation of $\tH_1$, $\tH_6$ and the fact that $H_6=i\pi/2$ implies that $\tH_5$ is oriented from $e^{L/2}$ to $-e^{L/2}$. So the orientation on $\tH_4$ is from $z$ to ${\bar z}$ where $z$ has positive imaginary part. But this implies that $H_5$ is purely real. So $\sigma_2=0$ which implies that $\sigma_1=-1$. 

As in the previous case, we test the above when $(a,b,c,d)=(1,0,1,1)$. This time hexagon $\sH$ looks as in figure \ref{fig:perp} below. In particular, the real part of $\tH_5$ is positive. With the formula for $H_5$ from theorem \ref{thm:main} this implies the lemma.

\begin{figure}[h]
 \begin{center}
 \ \psfig{file=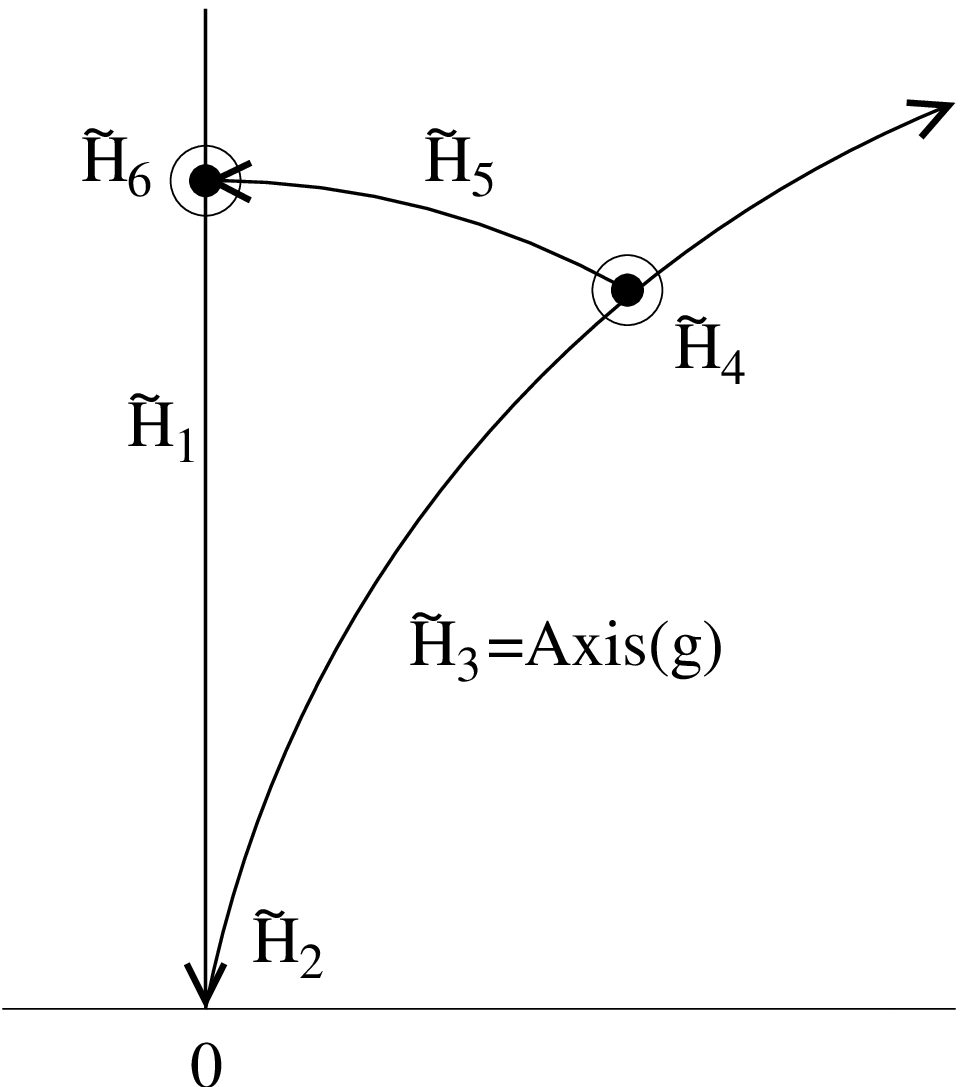, height=2.8in, width=2.3in}
 \caption{The hexagon $\hH$ in the case $M=-i\pi/2$ and $J=-ie^{L/2}$.}
 \label{fig:perp}
 \end{center}
 \end{figure}

\end{proof}


\begin{thebibliography}{99}

%\bibitem[Benedetti, Petronio]{BeP} R. Benedetti and C. Petronio, \textit{Lectures on Hyperbolic Geometry}, Springer-Verlag, Berlin, 1992.

%\bibitem[Agol, Long, Reid]{ALR} I. Agol, D.D. Long, A. Reid \textit{The Bianchi groups are separable on geometrically finite subgroups}, Annals of Math.\textbf{153} (2001), 599-621.

%\bibitem[Bekka \& Mayer]{Bek} M. Bekka, M. Mayer, \textit{Ergodic Theory and Topological Dynamics of Group Actions on Homogeneous Spaces}, London Math. Soc. Lecture Notes, 269 (2000), Cambridge University Press.

%\bibitem[Bestvina]{Bes} M. Bestvina, \textit{Questions in Geometric Group Theory}, http://www.math.utah.edu/~bestvina

%\bibitem[Bowen, L]{LBow} L. Bowen, \textit{Periodicity and Circle Packings
%  of the Hyperbolic Plane}, Geom. Dedicata. 102 (2003). pp213-236.

\bibitem[Bowen, R]{RBow} R. Bowen, \textit{The equidistribution of closed geodesics}. Amer. J. Math. 94 (1972), 413--423.  

\bibitem[Buser]{Bus} P. Buser, \textit{Geometry and Spectra of
  Riemann Surfaces}. Progress in Mathematics, 106. Birkhauser Boston
  Inc., Boston, MA, 1992. xiv + 454pp.

\bibitem[Ehrenpreis]{Ehr} L. Ehrenpreis, \textit{Cohomology with
  Bounds}. 1970. Symposia Mathematica. Vol IV (INDAM Rome 1968-69). pp389-395. 

\bibitem[Gendron]{Gen} T. Gendron, \textit{The Ehrenpreis Conjecture
  and the Moduli-Rigidity Gap}, Complex Manifolds and Hyperbolic
  Geometry (Guanojuato 2001), 207-229.

%\bibitem[Gordon, Long, Reid]{GLR} C. McA. Gordon, D. D. Long, A. Reid. \textit{Surface subgroups of Coxeter and Artin groups}, J. Pure and Applied Algebra 189 (2004), 135-148.

\bibitem[Hedlund1]{Hed} G.A. Hedlund, \textit{Fuchsian groups and Transitive Horocycles}, Duke Math J., (2), 530-543, 1936.

%\bibitem[Hedlund2]{Hed2} G.A. Hedlund, \textit{Dynamics of Geodesic Flows}, Bull. Am. Math. Soc. 40, (1939), 241-260. 

\bibitem[Fenchel]{Fen} W. Fenchel, \textit{Elementary Geometry in Hyperbolic Space}, Walter de Gruyter, Berlin, 1989.

%\bibitem[Li]{Li} Li, Tao Immersed essential surfaces in hyperbolic 3-manifolds.  Comm. Anal. Geom.  10  (2002),  no. 2, 275--290.

\bibitem[Ratcliffe]{Rat} J. Ratcliffe, \textit{Foundations of Hyperbolic
 Manifolds}, Springer-Verlag, New York, 1994.

\bibitem[Ratner]{Ra5} M. Ratner, \textit{Raghunathan's topological conjecture and distributions of unipotent flows}, Duke J. Math. (63), 235-280, 1991.

% \bibitem[Thurston]{Thu} W. Thurston, \textit{The Geometry and Topology of 3-manifolds}. Lecture Notes, Princeton.

%\bibitem[Waldhausen]{Wal} F. Waldhausen, \textit{The Word Problem in
%  Fundamental Groups of Sufficiently Large $3$-Manifolds}. Ann. of
%  Math. 88 (1968), pp272-280.

\bibitem[Zimmer]{Zim} R. Zimmer \textit{Ergodic Theory and Semisimple Groups}. Monographs in Mathematics, 81. Birkhäuser Verlag, Basel, 1984. x+209 pp. 


\end{thebibliography}
\end{document}